\documentclass[11pt]{amsart}
\usepackage[T1]{fontenc}
\usepackage[utf8]{inputenc}
\usepackage{import}
\usepackage{verbatim}
\usepackage[colorlinks=true,citecolor=blue,linkcolor=blue]{hyperref}
\usepackage{a4wide}
\usepackage{mathpazo}
\usepackage{epigraph}
\usepackage{xcolor}
\usepackage{mathrsfs}
\usepackage{amssymb, amsmath}
\usepackage{enumitem}
\usepackage{faktor}
\usepackage{tikz}
\newcommand{\Z}{{\mathbf{Z}}}
\newcommand{\N}{{\mathbf{N}}}

\newcommand{\R}{{\mathbf{R}}}

\newcommand{\1}{{\textbf{1}}}

\newcommand{\GL}{\mathbf{GL}}

\newcommand{\Span}{\mathrm{span}}
\newcommand{\supp}{\mathrm{supp}}

\newcommand{\Aut}{\mathrm{Aut}}
\newcommand{\I}{\mathrm{I}}
\newcommand{\Id}{\mathrm{Id}}
\newcommand{\Homeo}{\mathrm{Homeo^+}}
\newcommand{\TV}{\mathrm{TV}}
\newcommand{\RUC}{\mathcal{C}_{ru}^b}
\newcommand{\LUC}{\mathcal{C}_{lu}^b}
\newcommand{\UC}{\mathcal{C}_{u}^b}
\newcommand{\CB}{\mathcal{C}^b}
\newcommand{\CU}{\mathcal{C}^b_u}
\newcommand{\M}{\mathcal{M}}
\newcommand{\overp}{\overline{p}}

\theoremstyle{plain}
\newtheorem{thm}{Theorem}
\newtheorem*{thm*}{Theorem}
\newtheorem{lem}[thm]{Lemma}
\newtheorem{prop}[thm]{Proposition}
\newtheorem{cor}[thm]{Corollary}
\theoremstyle{definition}
\newtheorem*{defn*}{Definition}
\newtheorem{defn}[thm]{Definition}
\newtheorem{example}[thm]{Example}
\newtheorem*{example*}{Example}
\newtheorem{rem}[thm]{Remark}

\newtheorem{scholium}[thm]{Scholium}

\newtheorem{nota}[thm]{Notation}
\newtheorem*{rem*}{Remark}
\begin{document}

\title{Invariant Integrals on Topological Groups}

\begin{abstract}
We generalize the fixed-point property for discrete groups acting on convex cones given by Monod in \cite{monod} to topological groups. At first, we focus on describing this fixed-point property from a functional point of view, and then we look at the class of groups that have it. Finally, we go through some applications of this fixed-point property. To accomplish these tasks, we introduce a new class of normed Riesz spaces that depend on group representation.
\end{abstract}
\author{Vasco Schiavo}

\address{EPFL, Switzerland}

\date{Lausanne, October 2021}

%

\maketitle

\setcounter{tocdepth}{1}
\tableofcontents


\section{Introduction and Results}


In \cite{monod}, Monod defined a fixed-point property for discrete (abstract) groups acting on convex cones in locally convex vector spaces. He gave different characterizations of this fixed-point property, and he studied the class of groups that enjoy it. He pointed out that not all discrete groups have this fixed-point property, but those with it display interesting functional and hereditary properties. One of the most interesting features of Monod's fixed-point property is that it is strictly related to the concept of amenability. In fact, a discrete group that has this fixed-point property is supramenable and hence amenable. 

\medskip
Monod's fixed-point property was studied only for discrete groups. Therefore, a natural question is: is it generalizable to topological groups?

\medskip
There are essentially two factors that have led us to study this problem:
    \begin{itemize}
        \item[(I)] It is possible to define a good notion of amenability for topological groups.
        \item[(II)] In 1976, Jenkins defined in \cite{jenkinsexponential} a fixed-point property (called \textit{property F}) for locally compact groups while studying locally compact groups of subexponential growth. It turned out that Monod's fixed-point property and the one of Jenkins are equivalent when considering discrete groups, see \cite[Subsection 10.C]{monod}.
    \end{itemize}

\medskip
Therefore, taking inspiration from the theory of amenable topological groups and trying to unify the works of Jenkins and Monod, we define a generalization to topological groups of Monod's fixed-point property. After that, we start investigating it, trying to understand what difficulties and differences arise when going from the discrete to the topological case. First, we focus on describing this fixed-point property in a functional framework, i.e., using invariant functionals,  and then we study the class of groups which has it. Specifically, we try to understand this fixed-point property's hereditary and stability properties and what kind of groups are members of this class. Finally, we go through some applications of this fixed-point property.

\subsection{A topological fixed-point} The setting to define this fixed-point property is a locally convex vector space $E$, a non-empty weakly complete proper convex cone $C\subset E$ and a representation $\pi$ of $G$ on $E$ that leaves $C$ invariant.

\medskip
When we say that $G$ has a representation on a non-empty cone $C$ in a locally convex vector space $E$, we mean that $G$ has a linear representation on $E$, which leaves $C$ invariant.

\medskip
However, we have to put more conditions on the representation to avoid only finite groups enjoying this fixed-point property. In fact, even the action of $\R^*$, the multiplicative group of the real numbers, on the vector space $\R$ by scalar multiplication preserves the cone $C = [0,+\infty)$ but it has no non-zero fixed-point, since every non-zero orbit is unbounded. 

\medskip
At the same time, we would like to include the topology of the group in the data of the fixed-point property. 

\medskip
To satisfy both requests, we ask for a non-zero element  $x_0 \in C$ such that the orbital action of $G$ on $x_0$ is bounded right-uniformly continuous, i.e., for every $U \subset E$ neighborhood of the origin there is $V\subset G$ a neighborhood of the identity such that for every $a\in V$: $\pi(ag)x_0-\pi(g)x_0\in U$ for every $g\in G$ and such that the orbit of $x_0$ is a bounded subset of $E$. Recall that a subset of a locally convex vector space is said \textbf{bounded} if it is absorbed by every neighborhood of the origin. If there is a non-zero element $x_0\in C$ such that the orbital action of $G$ on $x_0$ is bounded right-uniformly continuous, then we say that the representation of $G$ on $C$ is \textbf{locally bounded right-uniformly continuous}. Note that the choice of such continuity condition is inspired by the Day-Rickert fixed-point criteria for amenability, see \cite[Theorem 4.2]{rickert}. 

\medskip
However, it was showed that any infinite (discrete) group acts on some locally compact space without fixing any non-zero Radon measure (\cite[Proposition 4.3]{matuirordam}). In contrast, every cocompact action of a (discrete) supramenable group on a locally compact space fixes a non-zero Radon measure (\cite[Proposition 2.7]{superamenablegroups}). The cocompactness of the action can be characterized by a property of the group's representation on the cone of Radon measures on the locally compact space in the following way. We say that a representation $\pi$ is of \textbf{cobounded} type if there is an element $\lambda \in E'$ in the topological dual of $E$ which is $G$-dominating, i.e., for every $\lambda'\in E'$ there are $g_1,...,g_n\in G$ such that $\pm \lambda ' \leq \sum_{j=1}^n \pi^*(g_j)\lambda$, where $\pi^*$ is the adjoint of the representation $\pi$ and the vector ordering on $E'$ is the one induced by the cone $C$. We recall that the order on $E'$ induced by $C$ is given by the relation: $\lambda \in E'$ is positive if and only if $\lambda(c) \geq 0$ for every $c\in C$.

\medskip
Finally, we propose the following generalization of Monod's fixed-point property for cones: 

\begin{defn}
Let $G$ be a topological group. We say that $G$ has \textbf{the fixed-point property for cones} if every representation of $G$ on a non-empty weakly complete proper convex cone $C$ in a (Hausdorff) locally convex vector space $E$ which is locally bounded right-uniformly continuous and of cobounded type has a non-zero fixed-point.
\end{defn}

\begin{rem}
Monod already proposed in \cite[Example 38]{monod} a possible generalization of his fixed-point property. However, the one he gave is not equal to ours in the non-locally compact case. We will motivate, and hopefully, we will persuade, the readers that our fixed-point property is the \textit{right one} to generalize the work of Monod.
\end{rem}

\subsection{Invariant functionals} It is customary and advantageous to translate a fixed-point property for groups into a functional analytic framework. Therefore, we also want to describe the fixed-point property for cones in this way. The functional which encodes its features is called an invariant normalized integral.

\medskip
Let $\RUC(G)$ be the space of all bounded right-uniformly continuous functions and let $G$ act on it by left-translation, i.e., $gf(x) = f(g^{-1}x)$ for every $g,x\in G$ and $f\in \RUC(G)$. Fix a non-zero positve function $f\in \RUC(G)$ and define the subspace $\RUC(G,f)$ of all functions in $\RUC(G)$ which are $G$-dominated by a finite number of translates of $f$. An invariant normalized integral on $\RUC(G,f)$ is a linear map $\I$ defined on $\RUC(G,f)$ which is positive, $G$-invariant and such that $\I(f)=1$. If the space $\RUC(G,f)$ admits an invariant normalized integral for every non-zero positive $f\in \RUC(G)$, then we say that \textbf{$G$ has the invariant normalized integral property for $\RUC(G)$}.

\begin{thm}\label{FPC equivalente integral}
Let $G$ be a topological group. Then the following are equivalent:
    \begin{itemize}
        \item[a)] the group $G$ has the fixed-point property for cones;
        \item[b)] the group $G$ has the invariant normalized integral property for $\RUC(G).$
    \end{itemize}
\end{thm}

The definition of the invariant normalized integral property is easily adapted to general function spaces. Therefore, it is natural to ask when the fixed-point property for cones is equivalent to the invariant normalized integral property on the classical Banach spaces $\CB(G)$, $\LUC(G)$, $\RUC(G)$ and $\CU(G)$. Note that a similar question was already studied and solved in the theory of invariant means (cf. \cite[Section $2$]{greenleaf}).

\medskip
We are far from an answer for the general topological case. However, we were capable to give a complete answer for locally compact groups.

\begin{thm}\label{Equivalence of integral on different spaces}
Let $G$ be a locally compact group. If $G$ has the invariant normalized integral property for one of the following Banach spaces
\begin{align*}
    L^{\infty}(G),\; \CB(G),\; \LUC(G),\; \RUC(G) \text{ or }\CU(G),
\end{align*}
then $G$ has the invariant normalized integral property for all the others. 
\end{thm}

In particular, this last theorem generalizes the main theorem of amenable locally compact groups theory. Namely,  \cite[Theorem $2.2.1$]{greenleaf}. This result states the equivalence of the existence of invariant means for the above classical Banach spaces when considering amenable groups.

\subsection{Looking at smaller spaces} Sometimes, it is not so pleasant to search for invariant normalized integrals on spaces of the form $\RUC(G,f)$. Can we instead focus on invariant functionals on smaller subspaces? A question already asked by Greenleaf in his famous monograph \cite[p.18]{greenleaf} is: given a non-zero $f\in \RUC(G)_+$, does the existence of a non-zero invariant functional on the subspace $\Span_\R(Gf)$ imply the existence of an invariant normalized integral on the space $\RUC(G,f)$? The essential information to ensure the existence of a non-zero invariant functional on $\Span_ \R(Gf)$ is contained in the following property of $f$.

\begin{defn}
  We say that a topological group $G$ has \textbf{the translate property for} $\RUC(G)$ if for every non-zero  $f \in \RUC(G)_+$ whenever
    \begin{align*}
        \sum_{j=1}^nt_jg_jf \geq 0 \quad \implies \quad \sum_{j=1}^nt_j \geq 0,
    \end{align*}
 where $t_1,...,t_n\in \R$ and $g_1,...,g_n\in G$.
\end{defn}

One of the first who studied this question was Rosenblatt, which was able to give a first partial answer for the discrete case in his PhD thesis (\cite[Corollary 1.3]{rosenblatt}). A complete answer for the discrete case was finally given by Monod (\cite[Corollary 19]{monod} and \cite[Corollary 20]{monod}). We were capable of giving a partial answer for the topological case (Proposition \ref{prop TP implies IIP for top groups}) and a complete one for locally compact groups. Indeed,

\begin{thm}\label{theorem TP imples measurably integral}
    Let $G$ be a locally compact group. If $G$ has the translate property for $\RUC(G)$, then $G$ has the invariant normalized integral property for $\RUC(G)$.
\end{thm}

Moreover, in the locally compact case, we answered Greenleaf's question not only for the space $\RUC(G)$ but for all the classical ones.

\begin{thm}\label{TP on E equivalent to NII on E}
Let $G$ be a locally compact group. Then the translate property and the invariant normalized integral property are equivalent for any of the following Banach spaces:
\begin{align*}
    L^{\infty}(G),\; \CB(G),\; \LUC(G),\; \RUC(G) \text{ and }\CU(G).
\end{align*}
\end{thm}

Thanks to these two last results, it was also possible to show the equivalence of the translate property on these spaces. 

\begin{thm}\label{Equivalence of TP on different spaces}
Let $G$ be a locally compact group. If $G$ has the translate property for one of the following Banach spaces 
\begin{align*}
    L^{\infty}(G),\; \CB(G),\; \LUC(G),\; \RUC(G) \text{ or }\CU(G),
\end{align*}
then $G$ has the translate property for all the others.
\end{thm}

\subsection{Groups with the fixed-point property for cones} After studying the functional properties of the fixed-point property for cones, we focus on studying the class of groups that enjoy it. This family is a subclass of the family of amenable groups. Unfortunately, it behaves not as good as it. This is because amenable groups are characterized by a fixed-point property on convex compact sets, which are easier to handle than cones. Nevertheless, it was possible to show some closure properties.

\begin{thm}\label{thm calss of top group with FPC}
    The class of topological groups with the fixed-point property for cones is closed by:
        \begin{itemize}
            \item[a)] taking open subgroups;
            \item[b)] surjective continuous homomorphism;
            \item[c)] taking extensions by finite groups.
        \end{itemize}
\end{thm}

It is immediately striking that the statement \textit{"being closed under extension"} is not included in the previous theorem. This is because it is not true, even for discrete groups (see Subsection \ref{Obstruction to the fixed point property on convex cones}). Moreover, it stands out that the assertion of being closed under directed union and taking the closure is not present. We were only able to prove some weaker results in the topological case (Propositions \ref{dense subgroup of topological with fixed-point property} and \ref{directed unin of topological subgroup}). Nevertheless, we have positive results for locally compact groups. 

\begin{thm}\label{ thm in loc cpt for dense/union/lattices}
Let $G$ be a locally compact group. Then:
    \begin{itemize}
        \item[a)] if $G$ is the directed union of closed subgroups with the fixed-point property for cones, then $G$ has the fixed-point property for cones;
        \item[b)] if $G$ has a dense subgroup with the fixed-point property for cones, then $G$ has the fixed-point property for cones.
    \end{itemize}
\end{thm}

We want to stress that the dense subgroup of point b) is not necessarily locally compact. Therefore, we should be careful about which characterization of the fixed-point property for cones we use for the ambient group and the dense subgroup. 

\medskip
Moreover, the following holds for locally compact groups.

\begin{thm}\label{Closed group have the FPC}
Suppose that $G$ is a locally compact group with the fixed-point property for cones and let $H\leq G$ be a closed subgroup. Then $H$ has the fixed-point property for cones. 
\end{thm}

We already mentioned that general group extensions do not preserve the fixed-point property for cones. Despite that, there are interesting extensions that preserve it for locally compact groups. 
Results in this directions were already given by Monod in the discrete case (see \cite[Theorem $8$ - (3), (4) \& (5)]{monod}).

\begin{thm}\label{cartesian product with subexponential growth}
Let $G$ be a locally compact group. Then:
    \begin{itemize}
        \item[a)] if $G$ is the direct product of a group with the fixed-point property for cones and one of subexponential growth, then $G$ has the fixed-point property for cones;
        \item[b)] if $G$ is the extension of a compact group by a group with the fixed-point property for cones, then $G$ has the fixed-point property for cones.
    \end{itemize}
\end{thm}

In particular, point a) of the above theorem shows that every locally compact group with subexponential growth has the fixed-point property for cones. Whether the converse is true is an open problem, even for discrete groups.

\subsection{Organization} In Section \ref{Chapter Ddominating and majorizing}, we study dominated normed Riesz spaces and positive invariant functionals on them. This part represents the toolbox for developing the theory in the next sections. 

In section \ref{from discrete to topological}, we investigate the fixed-point property for cones in the context of topological groups. We do the link with invariant integrals proving Theorem \ref{FPC equivalente integral}, and then explore Greenleaf's question in the topological case.

Section \ref{Chapter the locally compact case} is dedicated to locally compact groups. The demonstration of Theorem \ref{Equivalence of integral on different spaces} is provided here. Moreover, a complete answer to Greenleaf's question in the locally compact case is given with the proofs of Theorems \ref{theorem TP imples measurably integral} and \ref{TP on E equivalent to NII on E}.

In the short Section \ref{section equivalent fixed-point properties}, we use the results of Section \ref{Chapter the locally compact case} to get rid of some technical details when working with the fixed-point property for cones for locally compact groups. Those results will be helpful later to study the class of locally compact groups which enjoy the fixed-point property for cones.

The goal of section \ref{section hereditary properties}  is to study the class of groups which have the fixed-point property for cones and to understand which hereditary properties this class has. The proofs of the Theorems \ref{thm calss of top group with FPC}, \ref{ thm in loc cpt for dense/union/lattices}, \ref{Closed group have the FPC} and \ref{cartesian product with subexponential growth} are given here. Moreover, examples of groups with, and without, the fixed-point property are presented.

In section \ref{section fixing radon measures}, we study the relationship between non-zero invariant Radon measures and the fixed-point property for cones giving applications in problems where an invariant Radon measure is required. 

\subsection{Let's fix some notation} The capital letter $G$ always means a group with some specified topology. An abstract group is nothing but a group endowed with the discrete topology.

\medskip
Every vector space is to be considered real. The notation $E'$ and $E^*$ are used for the topological dual, the set of all continuous linear functionals on $E$ with respect to some given topology, and the algebraic dual, the set of all linear functionals on $E$, respectively.

\medskip
Every locally compact topology and every locally convex topology on a vector space is Hausdorff.

\bigskip\noindent
\textbf{Acknowledgements.} I'm thankful to my advisor N. Monod for introducing me to this topic and helping me with his advice and knowledge. Moreover, I'm grateful to M. Gheysens for many enlightening conversations.


\section{Dominating Riesz Spaces}\label{Chapter Ddominating and majorizing}

The purpose of this section is to recall some notions of the theory of ordered vector spaces and to develop new tools to work with a particular class of normed Riesz spaces, namely the dominated ones. Standard references for ordered vector spaces are \cite{banachlattice}, \cite{banachlatticeandposoperator} and the encyclopedic book \cite{infinitedim}. Regarding convex cones in vector spaces, we refer to the book \cite{conesandduality}.

\medskip
In this section, $G$ is always a topological group. 

\subsection{Order \& vectors}
An \textbf{ordered vector space} $(E, \leq)$ is nothing but a vector space equipped with a partial order $\leq$ such that for every $x,y,z\in E$ and every $c \in \R_+$ the relation $x \leq y$ implies $x+z \leq y+z $ and $cx\leq cy$. We only write $E$ instead of $(E, \leq)$, since it will always be clear which order we are considering. The set of \textbf{positive vectors} $E_+$ of $E$ is the set of all $x\in E$ such that $x \geq 0$. A linear functional $T$ between two ordered vector spaces $E$ and $V$ is said a \textbf{positive functional} if it sends positive vectors to positive vectors, i.e., if $v\in E_+$, then $T(v)\in V_+$. The functional $T$ is said \textbf{strictly positive} if it sends non-zero positive vectors to non-zero positive vectors, i.e., if $v\in E$ such that $v > 0$, then $T(v)> 0$. The set of all positive vectors of an ordered vector space $E$ is also called the \textbf{positive cone} of $E$.

\medskip
A \textbf{cone} $C$ is a non-empty subset of $E$ which is additive and positive homogeneous. The cone $C$ is said \textbf{proper} if $C\cap(-C) = \{ 0 \}$. If $C$ is a proper cone in a vector space $E$, then the binary relation on $E$ defined by
    \begin{align*}
        x \leq y \iff y-x\in C
    \end{align*}
makes $E$ an ordered vector space with positive cone $E_+ = C$. In fact, there is a one-to-one correspondence between proper cones and ordering on vector spaces (\cite[§1.1]{conesandduality}).

\medskip
An action of $G$ on an ordered vector space $E$ is always a representation of $G$ into the set of positive linear automorphisms of $E$, i.e., a group homomorphism $\pi : G \longrightarrow \Aut_L(E)$ such that whenever $v\leq w$ then $\pi(g)v \leq \pi(g)w$ for every $g\in G$ and every $v,w \in E$. Here, $\Aut_L(E)$ denotes the set of all linear automorphisms of $E$. Asking that the representation of $G$ on $E$ is positive is equivalent to asking that the positive cone of $E$ is $G$-invariant. Generally, we write $gv$ rather than $\pi(g)v$ for the action of $g\in G$ on the vector $v\in E$.

\medskip
We say that an ordered vector space $E$ is a \textbf{Riesz space} if for every pair of vectors $v,w \in E$ their infimum $v \wedge w$, or equivalently their supremum $v\vee w$, exists in $E$. Consequently, on a Riesz space $E$ makes sense to define a notion of absolute value via the equation $|v|= (-v)\vee v$ for $v\in E$. A vector subspace $F$ of a Riesz space $E$ is called a \textbf{Riesz subspace} if for every $v,w\in F$ their infimum $v \wedge w$, or equivalently their supremum $v\vee w$, belongs to $F$. From the identity $v \vee w = \frac{1}{2}\big( v+w +|v-w| \big)$, see \cite[Thm. $1.17 (6)$]{conesandduality} for a proof, a vector subspace $F$ of a Riesz space $E$ is a Riesz subspace if and only if $v\in F$ implies that $|v|\in F$. A vector subspace $S$ of a Riesz space $E$ is called an \textbf{ideal} of $E$ if $|v|\leq |w|$ and $w\in S$ implies that $v\in S$. As an ideal is closed by taking absolute value, then every ideal is automatically a Riesz subspace

\medskip
A \textbf{normed Riesz space} is a Riesz space $E$ together with an order-preserving, or monotone, norm, i.e., a norm $||\cdot ||$ on $E$ such that whenever $|v|\leq |w|$ then $||v||\leq ||w||$ for every $v,w\in E$.

\medskip
Let $G$ be a topological group with a representation by positive linear isometries on a normed Riesz space $E$. Then there are two natural notions of continuity: \textbf{orbital continuity}, i.e., for every $v\in E$ the map $g\longmapsto gv$ is continuous, and \textbf{joint continuity}, i.e., the map $(g,v)\longmapsto gv$ is continuous w.r.t. the product topology on $G\times E$. However, the two notions coincide for a topological group acting by isometries on a normed space, see \cite[Corollary C.1.6]{max}. For this reason, we simply say that $G$ acts continuously on $E$ or that the representation of $G$ on $E$ is continuous. 

\subsection{Dominated spaces}

From now until the end of the section, we suppose that $(E, ||\cdot||)$ is a normed Riesz space and that $G$ has a representation $\pi$ on $E$ by positive linear isometries.

\medskip
A vector $d\in E$ \textbf{$G$-dominates} another vector $v\in E$, or $v$ is \textbf{$G$-dominated} by $d$, if there are $g_1,...,g_n\in G$ such that $|v|\leq \sum_{j=1}^n g_jd$.

\begin{defn}
For a non-zero vector $d\in E_+$, we define the set
    \begin{align*}
        (E,d) := \left\{ v \in E : |v| \leq \sum_{j=1}^ng_jd \;\text{ for some }\; g_1,...,g_n\in G \right\}.
    \end{align*}
In other words, $(E,d)$ is the set of all vectors of $E$ which are $G$-dominated by $d$.
\end{defn}

The vector $d$ is called the \textbf{$G$-dominating element} of $(E,d)$.

\begin{nota}
Let $\mathfrak{F}(X)$ be a function space on a set $X$ with a Riesz space structure, and let $f$ be a non-zero positive element of $\mathfrak{F}(X)$. Then we write $\mathfrak{F}(X,f)$ instead of $(\mathfrak{F}(X),f)$.
\end{nota}

\begin{prop}
The set $(E,d)$ is an ideal of $E$ for every non-zero $d \in E_+$. In particular, $(E,d)$ is a Riesz subspaces of $E$.
\end{prop}
\begin{proof}
We start showing that $(E,d)$ is a linear subspace of $E$. Let $v,w\in (E,d)$, then there are and $g_1,...,g_n,h_1,...,h_m \in G$ such that $|v|\leq \sum_{j=1}^n g_jd$ and $|w|\leq \sum_{i=1}^m h_id $. Now, for every $\lambda\in \R$,
    \begin{align*}
        |\lambda v+w|\leq |\lambda||v|+|w| \leq |\lambda|\sum_{j=1}^n g_jd + \sum_{i=1}^m h_id.
    \end{align*}
Therefore, $\lambda v+ w \in (E,d)$. We can conclude that $(E,d)$ is a vector subspace of $E$. Moreover, the definition of $(E,d)$ implies that it is closed by taking absolute value. Hence, $(E,d)$ is a Riesz subspace of $E$.
\end{proof}

\subsection{Dominating norms} In the previous subsection, we gave a normed Riesz space and defined specific dominated Riesz subspaces of it. Those subspaces could be seen as normed Riesz subspaces only by restricting the norm of the ambient space on it. Unfortunately, this is useless for our future purposes. However, at least a norm arises naturally from the definition of $G$-dominated vector space with the properties we need.

\medskip
As before, let $(E, ||\cdot||)$ be a normed Riesz space and $\pi$ a representation of $G$ on $E$ by positive linear isometries.

\begin{defn}
 For every positive vector $d\in E,$ we define the possible infinite value
    \begin{align*}
        p_d(v) = \inf \left\{ \sum_{j=1}^nt_j : |v|\leq \sum_{j=1}^nt_jg_jd \;\text{ for some }\;t_1,...,t_n\in \R_+ \;\text{ and }\; g_1,...,g_n\in G \right\},
    \end{align*}
where $v\in E$. 
\end{defn}

The first thing to point out is that $p_d(v)$ is finite if and only if $v\in (E,d)$. Indeed, we have that $(E,d) = \left\{ v \in E : p_d(v) < \infty \right\}.$ Therefore, we proceed to study $p_d$ when restricted to $(E,d)$.

\begin{prop}\label{prop proprieta $G$-dominated norm}
Let $d\in E_+$ be a non-zero vector. Then:
    \begin{itemize}
        \item[a)] the pair $\left((E,d), p_d \right)$ is a normed Riesz space on which $G$ acts by positive linear isometries;
        \item[b)] the inequality  $||v|| \leq p_d(v)||d||$ holds for every $v\in (E,d)$;
        \item[c)] let $v,w\in (E,d)$ such that $v$ is $G$-dominated by $w$ and $w\geq 0$, then 
                        \begin{align*}
                            p_d(v)\leq p_w(v)p_d(w).
                        \end{align*}
    \end{itemize}
\end{prop}
\begin{proof}
The point a) is direct by the definition of $p_d$. 

Let's show point b). Let $\epsilon >0$ and take $v\in (E,d)$. Then there are $g_1,...,g_n\in G$ and $t_1,...,t_n\in \R_+$ such that $|v|\leq \sum_{j=1}^nt_jg_jd$ and $\sum_{j=1}^nt_j \leq p_d(v)+\frac{\epsilon}{||d||}$. Because the norm $||\cdot||$ is order-preserving, the following inequality holds:
    \begin{align*}
        ||v|| \leq \sum_{j=1}^n t_j||g_jd|| \leq \sum_{j=1}^n t_j ||d|| \leq p_d(v)||d|| +\epsilon. 
    \end{align*}
We can conclude that $||v|| \leq p_d(v)||d||$ as $\epsilon$ was chosen arbitrarily.  

For the proof of point c). Let $\epsilon > 0$. Then there are $t_1,...,t_n,c_1,...,c_m\in \R_+$, $g_1,...,g_n,h_1,...,h_m \in G$ such that $|v|\leq \sum_{j=1}^nt_jg_jw$ with $\sum_{j=1}^nt_j \leq p_w(v) + \sqrt{\epsilon}$ and $|w|\leq \sum_{k=1}^mc_kh_kd$ with $\sum_{k=1}^mc_k \leq p_d(w) + \sqrt{\epsilon} $. Therefore,
            \begin{align*}
                |v|   & \leq \sum_{j=1}^nt_jg_jw
                        \leq\sum_{j=1}^nt_jg_j \sum_{k=1}^m c_k h_k d
                        = \sum_{j=1}^nt_j\sum_{k=1}^m c_k g_jh_k d.
            \end{align*}
Since $p_d$ is also order-preserving, the following holds:
            \begin{align*}
                p_d(v) \leq \sum_{j=1}^nt_j\sum_{k=1}^m c_k \leq p_w(v)p_d(w)+ \sqrt{\epsilon}( p_w(v)+p_d(w)) + \epsilon.
            \end{align*}
We can conclude that $ p_d(v)\leq p_w(v)p_d(w)$ as $\epsilon$ was chosen arbitrarily.
\end{proof}

\begin{scholium}
Defining such type of dominating norms is a natural thing to do. In fact, there are a couple of examples where something similar was done. One is the functional proof of the existence of a Haar measure for every locally compact group, see for example \cite[Chap. VII]{bourbakihaar}. Here, dominating norms have been used on the set of compactly supported functions but without thinking of them as norms but as sub-additive maps with nice properties. Another \textit{quasi} example can be found in \cite[Theorem $2.55$]{conesandduality}, where we can see that $p_d$ is in some sense a Minkowsky functional for $G$-dominated Riesz spaces. 
\end{scholium}

\subsection{Invariant normalized integrals}
In this section, we look at positive functionals defined on dominated spaces. Such maps will follow us for the rest of the paper as they are the glue between the different sections. In particular, we will use those functionals in Section \ref{from discrete to topological} to describe the fixed-point property for cones in a functional framework. 

\medskip
As before, $(E, ||\cdot||)$ is a normed Riesz space and $\pi$ a representation of $G$ on $E$ by positive linear isometries.

\begin{defn}
Let $d\in E_+$ be a non-zero vector. We say that the Riesz space $(E,d)$ admits an \textbf{invariant normalized integral} if there exists a positive $G$-invariant and normalized functional on it, i.e., there is a functional $\I$ defined on $(E,d)$ such that $\I \geq 0$, $\I(d)=1$ and $\I(\pi(g)v) = \I(v)$ for every $g\in G$ and $v\in (E,d)$. 
\end{defn}

\begin{example}\label{example invariant normalized integral property}
    \begin{itemize}
        \item[1)] Every amenable topological group $G$ admits an invariant normalized integral, or a so-called mean, on the space $\RUC(G) = \RUC(G, \1_G).$
        
        \item[2)] First, note that $\mathcal{C}_{00}(G,\phi) = \mathcal{C}_{00}(G)$ for every non-zero $\phi \in \mathcal{C}_{00}(G)_+$. Indeed, if $\tilde{\phi} \in \mathcal{C}_{00}(G)$, then there are $g_1,...,g_n\in G$ such that $\1_{\supp(\tilde{\phi})} \leq \sum_{j=1}^n g_j \1_{\supp(\phi)}$. Therefore, $|\tilde{\phi}|\leq ||\tilde{\phi}||_\infty \sum_{j=1}^n g_j\phi$. Now consider a Haar measure $m_G$ for $G$. Then $m_G$ can be viewed as a strictly positive invariant functional on the Riesz space of compactly supported continuous functions $\mathcal{C}_{00}(G)$ (\cite[Chap. VII]{bourbakihaar}). Therefore, for a non-zero $\phi \in \mathcal{C}_{00}(G)_+$, the functional $\I = \frac{1}{m_G(\phi)} m_G$ is an invariant normalized integral on $\mathcal{C}_{00}(G,\phi) = \mathcal{C}_{00}(G)$. We can conclude that every locally compact group $G$ has the invariant normalized integral property for $\mathcal{C}_{00}(G)$.
    \end{itemize}
\end{example}

The inspiration to define $p_d$-norms comes from the following corollary. In fact, such norms are the good ones to generate a topology for which invariant normalized integrals are continuous.

\begin{cor}\label{continuity of the normalized invariant integral}
Let $E$ be a normed Riesz space, and let $d\in E_+$ be a non-zero vector. Then every invariant normalized integral $\I$ on $(E,d)$ is continuous for the norm $p_d$ and has operator norm equal to 1.
\end{cor}

\begin{proof}
Let $v\in (E,d)$ and $\varepsilon > 0.$ Then there are $t_1,...,t_n\in \R_+$ and $g_1,...,g_n\in G$ such that $|v| \leq \sum_{j=1}^n t_jg_jd$ and $\sum_{i=1}^n t_j \leq p_d(v)+ \frac{\varepsilon}{M}$. Thus,
    \begin{align*}
        \left| \I(v)\right| & \leq  \I\left(\sum_{i=1}^n t_jg_jd \right)
                              \leq \sum_{j=1}^n t_j \I(g_jd) =  \sum_{j=1}^n t_j \leq p_d(v) + \varepsilon.
    \end{align*}
Therefore, $||\I||_{op} \leq 1$ as $\varepsilon$ was chosen arbitrarily. Meanwhile, $\I(d) =1$. Hence, $||\I||_{op}=1$.
\end{proof}

We are ready to define the notion that will permit us to build the link between the fixed-point property for cones and invariant normalized integrals.

\begin{defn}
Let $E$ be a normed Riesz space and suppose that $G$ has a representation $\pi$ on $E$ by positive linear isometries. We say that \textbf{$G$ has the invariant normalized integral property for $E$} if there is an invariant normalized integral on $(E,d)$ for every non-zero $d\in E_+$.
\end{defn}

\begin{example}
    \begin{itemize}
        \item[1)] Let $G$ be a locally compact group and consider the natural action of $G$ on $L^1(G)$. Then the formula $\I(f) = \int_G fdm_G$ defines a strictly positive invariant functional on $L^1(G)$. Hence, every locally compact group $G$ has the invariant normalized integral property for $L^1(G)$.
        
        \item[2)] Every compact group $G$ has the invariant normalized integral property for $L^p(G)$ for $1\leq p\leq \infty$. This is because the inclusion $L^p(G) \subset L^1(G)$ holds for every $1 \leq p\leq \infty$ (\cite[Corollary 13.3]{infinitedim}).
        
        \item[3)] Let $G$ be a topological group and consider the normed Riesz space of almost periodic functions $\mathcal{AP}(G)$ on $G$, i.e., the spaces of all continuous bounded functions which have relatively compact orbit with respect to the $||\cdot||_\infty$-norm. Then $\mathcal{AP}(G)$ admits a strictly positive invariant functional, see \cite[§3.1]{greenleaf}. Hence, every topological group $G$ has the invariant normalized integral property for $\mathcal{AP}(G)$.
    \end{itemize}
\end{example}

\section{From Discrete to Topological}\label{from discrete to topological}

The primary goal of this section is to characterize the fixed-point property for cones for topological groups using the invariant normalized integral property. This is done by presenting the proof of Theorem \ref{FPC equivalente integral}. After that, we try to answer Greenleaf's question in the topological case by looking at the relationship between invariant normalized integrals and the translate property.

\medskip
For this section, $G$ is a topological group. 

\subsection{Invariant normalized integrals} Before illustrating the proof of Theorem \ref{FPC equivalente integral}, we recall the definition of some function spaces. 

\medskip
Let $\CB(G)$ be the space of all bounded continuous functions on $G$. We define the left-translation representation $\pi_L$ of $G$ on $\CB(G)$ by $\pi_L(g)f(x) = f(g^{-1} x)$, where $g,x\in G$ and $f\in \CB(G)$, and the right-translation representation $\pi_R$ of $G$ on $\CB(G)$ by $\pi_R(g)f(x) = f(xg)$, where $g,x\in G$ and $f\in \CB(G)$. We write $\RUC(G)$ for the subspace of all bounded right-uniformly continuous functions on $G$, i.e., the space of all $f\in \CB(G)$ such that the orbital map $g \longmapsto \pi_L(g)f$ is continuous with respect to the supremum norm on $\CB(G)$. Similarly, we write $\LUC(G)$ for the subspace of all bounded left-uniformly continuous functions on $G$, i.e., the space of all $f\in \CB(G)$ such that the orbital map $g \longmapsto \pi_R(g)f$ is continuous with respect to the supremum norm on $\CB(G)$. The space $\CU(G)$ of all bounded uniformly continuous functions on $G$ is defined as the intersection of $\RUC(G)$ and $\LUC(G)$.

\medskip
Unless otherwise specified, we always consider the left-translation representation on all of the above function spaces, and we write $gf$ instead of $\pi_L(g)f$ and $f_g$ instead of $\pi_R(g)f$ for $g\in G$ and $f\in \CB(G)$.

\medskip
We are ready to demonstrate Theorem \ref{FPC equivalente integral}. Our proof follows the one given by Monod in \cite[Theorem $7$]{monod}. But first, we have to take care of the locally bounded right-uniformly continuous condition. This is done in the following lemma.

\begin{lem}\label{lemma induced G-repr is of locally RUC type}
Let $f\in \RUC(G)$ be a non-zero positive function and let $g\in G$ such that $f(g)\neq 0$. Write $ev_g$ for the evaluation map at point $g$. Then the orbital map $ \omega: x\longmapsto \omega(x) = x\cdot ev_g$ is bounded right-uniformly continuous with respect to the weak-* topology on $\RUC(G,f)^*$.
\end{lem}
\begin{proof}
We begin by showing that $\omega$ is bounded. Therefore, let $U\subset \mathcal{C}_{ru}^b(G,f)^*$ be an open neighborhood of the origin. We can suppose that $U$ is of the form
    \begin{align*}
        U = \big\{ \psi \in \mathcal{C}_{ru}^b(G,f)^*: |\psi(f_j)| < \delta \quad \text{ for } j=1,...,n\big\},
    \end{align*}
where $\delta >0$ and $f_1,...,f_n\in \mathcal{C}_{ru}^b(G)$, as a neighborhoods basis at the origin for the weak-* topology is given by such sets (\cite[Chap. II §6 No.2]{bourbakiespace}). We have to show that there is $t\in \R_+$ such that $im(\omega)= G \cdot ev_g\subset tU$. Set $t>\frac{\max_{j}||f_j||_\infty}{\delta}$ and let $ev_h\in im(\omega)$, where $h\in G$. Then
    \begin{align*}
        |ev_h(f_j)|=|f_j(h)|\leq \max_{1\leq j \leq n} ||f_j||_{\infty} < \delta t.
    \end{align*}
We can conclude that $ev_g\in tU$. It is left to show that $\omega$ is right-uniformly continuous. Let $U\subset \mathcal{C}_{ru}^b(G,f)^*$ be an open neighborhood of the origin of the same form as above. Let $\epsilon= \delta >0$ and take $V\subset G$ an open neighborhood of the identity such that for every $a\in V:$
    \begin{align*}
        |f_j(ag)-f_j(g)|<\delta\quad \text{ for all }g\in G\text{ and } j=1,...,n.
    \end{align*}
This is possible because the $f_j's$ are in $\mathcal{C}_{ru}^b(G)$. Therefore,
    \begin{align*}
        |\omega(ag)(f_j) -\omega(g)(f_j)|  = | \big(ev_{ag}-ev_{g}\big)(f_j)| 
             = |f_j(ag)-f_j(g)| < \epsilon = \delta
    \end{align*}
for all $g\in G$ and $j=1,...,n$. This implies that $\omega(ag)-\omega(g)\in U$ proving the right-uniform continuity.
\end{proof}

We are therefore ready to provide a proof of Theorem \ref{FPC equivalente integral}.

\begin{proof}[Proof of Theorem \ref{FPC equivalente integral}]
We start by proving that a) implies b). Thus, let $f\in \RUC(G)$ be a non-zero positive function. We want to show that there is an invariant normalized integral $\I$ on $\RUC(G,f)$. Set $E=\RUC(G,f)^* $. Then $E$ is a Hausdorff locally convex topological vector space when equipped with the weak-* topology, see \cite[I §6 No.2 Remarque 1)]{bourbakiespace}. Consider now the cone $C =\RUC(G,f)^*_+$ of positive functionals on $E$. This cone is invariant and proper because spanned by positive elements. Recall that closed subspaces of complete spaces are complete (\cite[II §3 No.4 Proposition 8]{bourbakitopgen}) and that algebraic duals are weak-* complete (\cite[II §6 No.7]{bourbakiespace}). Therefore, $C$ is weakly complete as it is closed in the complete space $E$. We claim that the representation of $G$ on $C$ is of locally bounded right-uniformly continuous type. Indeed, let $x\in G$ such that $f(x) \neq 0$. Then the map $g \longmapsto gev_x$ is bounded and right-uniformly continuous by Lemma \ref{lemma induced G-repr is of locally RUC type}. Finally, the action is of cobounded type, since $\left( E, \text{weak-*}\right)'= \RUC(G,f)$ by \cite[II §6 No.1 Remarque 1)]{bourbakiespace}. Therefore, there is a non-zero $G$-invariant element $\I\in C$ thanks to the fixed-point property for cones of $G$. It is clear that $\I(f)\neq 0$. Indeed, let $\phi \in E_+$ such that $\I(\phi) > 0$. Then $\phi \leq \sum_{j=1}^n g_jf$ for some $g_1,...,g_n\in G$. Therefore, 
    \begin{align*}
        0<\I(\phi) \leq \I \left( \; \sum_{j=1}^n g_jf \, \right) = n\I(f).
    \end{align*}
After a normalization, $\I$ becomes an invariant normalized integral for $E$.

Let's show the inverse implication, i.e., b) implies a). Suppose that $G$ has a locally bounded right-uniformly continuous representation on a non-empty weakly complete proper cone $C$ in a locally convex space $E$ which is of cobounded type. Fix $x_0 \in C$ a non-zero vector with bounded right-uniformly continuous $G$-orbit, and define the positive linear operator
    \begin{align*}
        T: E' \longrightarrow \RUC(G), \quad \lambda \longmapsto T(\lambda)
    \end{align*}
where $T(\lambda)(g) = \lambda(gx_0)$ for every $g\in G$. The operator $T$ is well-defined because it is a composition of uniformly continuous functions and $T(\lambda)$ is bounded for every $\lambda\in E'$ because the set $Gx_0 \subset E$ is bounded and continuous linear functionals map bounded sets to bounded sets (\cite[III §1 No.3 Corollaire 1]{bourbakiespace}). Moreover, $T$ is also equivariant. Now, let $d\in E'$ be the $G$-dominating element given by the cobounded condition. Then $T$ maps $E'$ into $\RUC(G, T(d))$. On this last space, there is an invariant normalized integral by hypothesis. We can now use the same strategy employed in the implication $(4) \implies (1)$ of \cite[Theorem 7, pp. 77-78]{monod} to ensure the existence of a non-zero fixed-point in $C$.
\end{proof}

\subsection{The translate Property} We explore Greenleaf's question for topological groups. In other words, we investigate the relation between the translate property and the invariant normalized integral property. 

\medskip
First of all, note that the translate property for $\RUC(G)$ carries enough information to ensure the existence of an invariant mean on $\RUC(G)$, i.e.,  the topological group $G$ is amenable. In fact,

\begin{thm}\label{translate property implies amenability}
 Suppose that $G$ has the translate property for $\RUC(G)$, then $G$ is amenable.
\end{thm}

The proof follows the strategy of \cite[Theorem $18$]{monod} but with a significant variation. In fact, the core of Monod's proof was a theorem of Moore, which characterizes amenability in the discrete setting. Sadly, this theorem does not hold for topological groups. Luckily, there is a generalisation of it which fit our needs perfectly.

\begin{proof}[Proof of Theorem \ref{translate property implies amenability}.]
By \cite[Theorem 3.2 (b)]{lau}, it suffices to show that, for every non-zero function $f\in \mathcal{C}_{ru}^b(G)_+$, there is a mean $m$ on $\mathcal{C}_{ru}^b(G)$ such that $m(gf) = m(f)$ for every $g\in G$. Thus, fix a non-zero $f\in \mathcal{C}_{ru}^b(G)_+$ and consider the Banach subspace
    \begin{align*}
        D = \overline{\Span_{\R}\big\{ f-gf:g\in G\big\}}^{||\cdot||_{\infty}} \subset \mathcal{C}_{ru}^b(G).
    \end{align*}
Denote $Q = \faktor{\RUC(G)}{D}$ the quotient space with quotient norm $||\cdot||_{Q}$. We claim that $||\1_Q||_{Q} = 1$. Suppose it is not the case, then there is $\epsilon > 0,t_1,...,t_n\in \R$ and $g_1,...,g_n\in G$ such that $v = \sum_{j=1}^nt_j(f-g_jf)$ satisfies $|| \1_G -v||_{\infty} \leq 1-\epsilon.$ This means that
    \begin{align*}
        (\1_{G}-v)(x)\leq (1-\epsilon)\1_G(x) \iff \epsilon\1_{G}(x)\leq v(x) \quad \text{ for every }x\in G.
    \end{align*}
Let now $M := ||f||_{\infty}$ and note that
    \begin{align*}
        \frac{\epsilon}{M}f(x)\leq \epsilon\1_{G}(x)\leq v(x) \iff v(x)-\frac{\epsilon}{M}f(x)\geq 0
    \end{align*}
for every $x\in G$. But this element is in contradiction with the fact that $G$ has the translate property. Indeed, the sum of its coefficients is equal to $-\frac{\epsilon}{M}<0.$ Therefore, $||\1_Q||_{\infty} = 1$. By the Hahn-Banach Theorem, there is a continuous linear functional $m_Q$ on $Q$ such that $m_Q(\1_Q)=1$. Now define $m$ as the lift of $m_Q$. Then $m$ is of norm one and positive. Moreover, $m(f)=m(gf)$ for every $g\in G$, since it vanishes on $D$.
\end{proof}

The next natural question is: what kind of functionals does the translate property allow us to construct on linear subspaces of $\RUC(G)$?

\begin{prop}\label{translate property implic che I(G) non vuoto}
The following assertions are equivalent: 
    \begin{itemize}
        \item[a)] the group $G$ has the translate property for $\RUC(G)$; 
        \item[b)] for every non-zero $f\in \mathcal{C}_{ru}^b(G)_+$ there is $\psi \in \mathcal{C}_{ru}^b(G,f)^*_+$ such that $\psi(gf)=\psi(f)=1$ for every $g \in G$.
    \end{itemize}
\end{prop}

\begin{proof}
We show that b) implies a). Let $f\in \RUC(G)_+$ be a non-zero function and let $t_1,...,t_n\in \R$ and $g_1,...,g_n\in G$ such that $\sum_{j=1}^nt_jg_jf \geq 0$. Let $\psi$ as in the hypothesis of b). Then 
    \begin{align*}
       0 \leq  \psi \left(\, \sum_{j=1}^n t_jg_jf \right) = \sum_{j=1}^n t_j \psi(g_jf) = \sum_{j=1}^n t_j.
    \end{align*}
Therefore, $G$ has the translate property for $\RUC(G)$.
So let's prove the other direction. Fix a non-zero $f\in \mathcal{C}_{ru}^b(G)_+$ and define the linear functional
    \begin{align*}
        \omega : \Span_\R \left\{ gf: g\in G\right\}  \longrightarrow \R, \quad
                \sum_{j=1}^nt_jg_jf  \longmapsto \sum_{j=1}^nt_j
    \end{align*}
Note that $\omega$ is well-defined thanks to the translate property. Moreover, $\omega$ is positive, $G$-invariant and $\omega(f)=1$. 
As $\Span_\R \left\{ gf: g\in G\right\}$ is a $G$-dominating subspace of $\RUC(G,f)$, it is possible to extend $\omega$ in a positive way to $\RUC(G,f)$ by Kantorovich Theorem (\cite[Corollary $1.5.9$]{banachlattice}). This extension is the functional $\psi$ asked by point b).  
\end{proof}

A direct application shows that:

\begin{cor}\label{corollary INI property implies TP}
Suppose that $G$ has the invariant normalized integral property for $\RUC(G)$, then $G$ has the translate property for $\RUC(G)$.
\end{cor}

We introduce the following definition to understand in which case the converse also holds.

\begin{defn}
 Let $f \in \RUC(G)_+$ be a non-zero function and consider  $\RUC(G,f)$ equipped with the $p_f$-norm. Define
    \begin{align*}
        \mathcal{I}_f(G) = \big\{ \psi \in \mathcal{C}_{ru}^b(G,f)'_+ : \psi(gf) = \psi(f)= 1 \text{ } \forall g\in G \big\} \subset \mathcal{C}_{ru}^b(G,f)'.
    \end{align*}
\end{defn}

An application of Banach-Alaoglu Theorem (\cite[Theorem 5.105]{infinitedim}) gives the following: 

\begin{prop}\label{I(G) e compatto}
The set $\mathcal{I}_f(G)$ is convex and compact with respect to the weak-* topology for every non-zero function $f\in \mathcal{C}_{ru}^b(G)_+$.
\end{prop}

Note that, for a general topological group $G$, it is possible to find a non-zero $f\in \RUC (G)_+$ such that the set $\mathcal{I}_f(G)$ is empty. Take as an example a non-amenable group. However, if $G$ has the translate property for $\RUC(G)$, then $\mathcal
{I}_f(G)$ is non-empty for every non-zero positive $f\in \RUC(G)_+$ by  Proposition \ref{translate property implic che I(G) non vuoto}.

\medskip
The next proposition is the most satisfying answer to the Greenleaf's question that we could give in the topological case.

\begin{prop}\label{prop TP implies IIP for top groups}
Suppose that $G$ has the translate property for $\RUC(G)$, and let $f\in \RUC(G)$ be a non-zero positive function. If that the action of $G$ on $\RUC(G,f)$ is continuous with respect to the $p_f$-norm, then $\RUC(G,f)$ admits an invariant normalized integral.
\end{prop}
\begin{proof}
Recall that $\mathcal{I}_f(G)$ is non-empty because $G$ has the translate property for $\RUC(G)$. Moreover, the adjoint action of $G$ on $\mathcal{I}_f(G)$ is orbitally continuous with respect to the weak-* topology, as the action of $G$ on $\RUC(G,f)$ is $p_f$-continuous. But now $G$ is amenable by Theorem \ref{translate property implies amenability}. Therefore, there is a $G$-fixed-point in $\mathcal{I}_f(G)$, which is an invariant normalized integral on $\RUC(G,f).$ 
\end{proof}

This last proposition generalizes old answers to the Greenleaf's question for discrete groups. For example, the results \cite[Proposition $2.14$]{kellerhals} and \cite[Theorem $1.3.2$]{greenleaf} are corollaries of it.

\medskip
Unfortunately, it was not possible to show that the action of $G$ on $\RUC(G,f)$ is $p_f$-continuous for every non-zero $f\in \RUC(G)_+$. Nevertheless, we could show it for a special set of functions. Precisely, for those which are support-dominating. Recall that the support of a real function $f:X \longrightarrow \R$, defined on an arbitrary set $X$, is the subset $\supp(f) = \left\{ x\in X : f(x) \neq 0 \right\}$.

\begin{defn}
We say that a non-zero positive function $f\in \RUC(G)$ is \textbf{support-dominating} if there are $g_1,...,g_n\in G$ such that $\1_{\supp(f)}\leq  \sum_{j=1}^n g_jf$. 
\end{defn}

Note that in general $\1_{\supp(f)} \notin \RUC(G,f)$ or $\RUC(G)$.

\begin{lem}\label{la azione é orbitalli continuous}
Let $f\in \mathcal{C}_{ru}^b(G)_+$ be a non-zero support-dominating function. Then the action of $G$ on $\mathcal{C}_{ru}^b(G,f)$ is $p_f$-continuous.
\end{lem}

\begin{proof}
Fix a $\phi \in \mathcal{C}_{ru}^b(G,f)$ and let $(g_\alpha)_{\alpha}\subset G$ be a net such that $\lim_\alpha g_{\alpha} = e$. We want to show that $\lim_{\alpha}p_f(g_\alpha\phi -\phi)=0$. Let $\epsilon > 0$ and note that the inequality
    \begin{align*}
            |g_{\alpha}\phi-\phi| \leq ||g_{\alpha}\phi-\phi||_{\infty}\left(g_{\alpha}\1_{\supp(\phi)} +\1_{\supp(\phi)}\right).
    \end{align*}
holds for every $\alpha$.
As $f$ is support-dominating, there are $g_1,...,g_n\in G$ such that $\1_{\supp(\phi)} \leq \sum_{j=1}^n g_jf$. Therefore,
    \begin{align*}
            |g_{\alpha}\phi-\phi| \leq ||g_{\alpha}\phi-\phi||_{\infty}\left(g_{\alpha}\sum_{j=1}^n g_jf +\sum_{j=1}^n g_jf \right)
     \end{align*}
for every $\alpha.$ Because the representations of $G$ on $\RUC(G)$ is continuous for the $||\cdot||_\infty$-norm, there is $\alpha_0$ such that $|| g_\alpha \phi - \phi ||_\infty < \frac{\epsilon}{2n}$ for every $\alpha \succeq \alpha_0$. We can conclude that 
    \begin{align*}
            p_f( |g_{\alpha}\phi-\phi| )    \leq ||g_{\alpha}\phi-\phi||_{\infty}p_f\left(g_{\alpha}\sum_{j=1}^n g_jf +\sum_{j=1}^n g_jf \right)  
            \leq 2n ||g_{\alpha}\phi-\phi||_{\infty} < \epsilon
    \end{align*}
for every $\alpha \succeq \alpha_0$.
\end{proof}

In general, convergence in $||\cdot||_\infty$-norm doesn't imply convergence in $p_f$-norm. However, the two convergences coincide when the dominating element is support-dominating and with some other conditions.

\begin{prop}\label{prop conv of net in sup norm implies convergence in p_f}
Let $f\in \mathcal{C}_{ru}^b(G)_+$ be a non-zero support-dominating function, and let $(\phi_\alpha)_\alpha$ by a net in $\mathcal{C}_{ru}^b(G,f)$ which converges in $||\cdot||_{\infty}$-norm to $\phi \in \mathcal{C}_{ru}^b(G,f)$ and which has decreasing supports, i.e., $supp(\phi_\alpha) \supset supp(\phi_{\alpha'})$ if $\alpha' \succeq \alpha$. Then $(\phi_\alpha)_\alpha$ converges to $\phi$ in $p_f$-norm.
\end{prop}

\begin{proof}
Let $\epsilon>0$ and fix an $\alpha_0$. As $f$ is a support-dominating function, there are $g_1,...,g_n, y_1,...,y_m\in G$ such that $\1_{\supp(\phi)} \leq \sum_{j=1}^n g_jf$ and $\1_{\supp(\phi_{\alpha_0})} \leq \sum_{i=1}^m y_if$. Now the inequality

    \begin{align*}
        | \phi - \phi_{\alpha} | \leq || \phi - \phi_{\alpha} ||_{\infty} \left( \1_{\supp(\phi)} + \1_{\supp(\phi_{\alpha})} \right)
            \leq || \phi - \phi_{\alpha} ||_{\infty} \left(\, \sum_{j=1}^n g_jf + \sum_{i=1}^m y_if \right)
    \end{align*}
holds for every $\alpha \succeq \alpha_0$. Taking $\alpha'$ such that $|| \phi - \phi_{\alpha} ||_{\infty} < \frac{\epsilon}{n+m}$ and $\alpha' \succeq \alpha_0$, it is possible to conclude that
    \begin{align*}
        p_f(|\phi_{\alpha}-\phi|) \leq || \phi - \phi_{\alpha} ||_{\infty} p_f \left(\, \sum_{j=1}^n g_jf + \sum_{i=1}^m y_if \right) < \epsilon
        \quad \text{ for every } \alpha \succeq \alpha'.
    \end{align*}
\end{proof}


\section{The Locally Compact Case}\label{Chapter the locally compact case}
For this section,  let $G$ be a locally compact group and $m_G$ a fixed left-invariant Haar measure of $G$. 

\medskip
Write $L^\infty(G)$ for the space of all bounded (Borel) measurable functions on $G$ modulo null sets, and $L^1(G)$ for the space of all (Borel) measurable functions on $G$ which are $m_G$-integrable. If it is not specified, the letter $E$ is free to be one of the following Banach spaces:
    \begin{align*}
        \left\{ L^{\infty}(G),\, \CB(G),\, \LUC(G),\, \RUC(G),\, \CU(G)\right\}.
    \end{align*}
Write $\mathcal{C}_{00}(G)$ for the Riesz space of all compactly supported continuous functions on $G$. 

\medskip
On every of the above function spaces,  $G$ acts via the left-translation representation.

\subsection{Measures and convolution}

Let $\mathcal{B}(G)$ be the Borel $\sigma$-algebra of the group $G$, i.e., the $\sigma$-algebra generated by the open sets of the topology of $G$. A measure on $G$ is called a \textbf{Borel measure} if it has domain equal to $\mathcal{B}(G)$. A Borel measure is said locally finite if it takes finite values on compact sets. A locally finite Borel measure which is both inner and outer regular, is called a \textbf{regular Borel measure}. A \textbf{signed measure} is a measure which can also take negative values, and a \textbf{finite measure} is one which only takes finite values.

Write $\M (G)$ for the vector space of \textbf{all signed finite regular Borel measures on $G$}. Then $\M(G)$ becomes a Riesz space when equipped with the order given by
    \begin{align*}
        \mu \leq \lambda \quad \iff \quad \mu(A) \leq \lambda(A) \text{ for all }A\in \mathcal{B}(G).
    \end{align*}
We reefer to \cite[p. $22$]{conesandduality} for details.

\medskip
The \textbf{total variation} of a measure $\mu \in \M (G)$ is given by $||\mu||_\TV= |\mu|(G)$. This formula defines a norm on $\M(G)$ which turns $\M (G)$ into a Banach space (\cite[Proposition $4.1.8$]{measurecohn}). In particular, $\M (G)$ is a normed Riesz space, since the total variation norm is order-preserving. 
Note that $\M(G)$ is generated by positive measures thanks to the Jordan Decomposition Theorem, which says that every measure in $\M (G)$ can be written as the difference of two positive measures (\cite[Corollary $4.1.6$]{measurecohn}).

\medskip
Let $\mu$ and $\upsilon$ belong to $\M (G)$ and define their \textbf{convolution} as 
    \begin{align*}
        (\mu*\upsilon) (A) = \int_G \upsilon (x^{-1} A)d\mu(x) = \int_G \mu (A y^{-1})d\upsilon (y),
    \end{align*}
where $A \in \mathcal{B}(G)$. The preceding expression is well-defined and belongs to $\M (G)$ thanks to the Jordan Decomposition Theorem and \cite[Lemma $9.4.5$]{measurecohn}. The Banach space $\left( \M (G), ||\cdot||_\TV \right)$ with the convolution as multiplication is a unital Banach algebra (\cite[Proposition $9.4.6$]{measurecohn}).

\medskip
Recall that the space $L^1(G)$ is isometrically isomorphic as a Banach algebra to the ideal \textbf{$\M (G)_a$ of all signed finite regular Borel measures which are absolutely continuous with respect to $m_G$} (\cite[Theorem $A.1.12$]{runde}). The isometric isomorphism is given by the map $f\mapsto \mu_f$, where the measure $\mu_f$ is defined as $\mu_f(A)= \int_A f(g) d\mu(g)$. Note that here ideals are considered in the algebraic sense.

\medskip
There is a natural left \textbf{convolution-action} of $\M (G)$ on $L^\infty(G)$ given by
    \begin{align*}
        (\mu * f)(g) = \int_G f(x^{-1}g)d\mu(x) \quad \text{ for $m_G$-almost every }g \in G,
    \end{align*}
where $\mu \in \M (G)$ and $f\in L^{\infty}(G).$ Then $\mu*f$ is a bounded measurable function. Indeed,
    \begin{align*}
        ||\mu * f ||_{\infty} \leq || \mu ||_\TV ||f||_\infty \quad\text{ for every }\mu \in \M (G) \text{ and every }f\in L^{\infty}(G).
    \end{align*}
See \cite[Theorem $20.12$]{abstractharmonic} for more details. 

\medskip
If $\lambda \in \M(G)_a$ and $\theta$ is its corresponding function in $L^1(G)$, then
    \begin{align*}
        (\lambda*f)(g) = \int_G f(x^{-1}g)d\lambda(x) = \int_G f(x^{-1}g)\theta(x)dm_G(x) = (\theta*f)(g)
    \end{align*}
for every $g\in G$ and $f\in L^\infty(G)$.

\begin{nota}\label{notation map m_f}
For a function $\theta \in L^1(G)$, we write $\lambda_\theta$ to refer to the measure in $\M(G)_a$ given by the equation $d\lambda_\theta = \theta  dm_G$. Conversely, for a measure $\lambda \in \M(G)_a$, we write $\theta_\lambda$ to refer to the preimage of $\lambda$ in $L^1(G)$.
\end{nota}

Before going further, recall that, for a locally compact group $G$, the following isomorphisms hold:
    \begin{align*}
        \RUC(G) & = \left\{ \phi * f : \phi \in L^1(G) \text{ and } f\in L^\infty(G) \right\}    \\
        \LUC(G) & = \left\{ f * \phi : f \in L^\infty(G) \text{ and } \phi \in L^1(G)\right\}   \\
        \CU(G)  & = \left\{ \phi_1 * f * \phi_2 : \phi_1, \phi_2 \in L^1(G) \text{ and } f \in L^\infty(G) \right\}.
    \end{align*}
This is essentially due to Cohen-Hewitt Factorization Theorem (\cite[(32.22)]{abstractharmonic2}). We refer to \cite[(32.45)]{abstractharmonic2} for a proof.

\begin{lem}\label{invariant for M(G)}
 The Banach space $E$ is $\M (G)$-invariant.
\end{lem}
\begin{proof}
The space $L^\infty (G)$ is $\M(G)$-invariant by definition of the convolution-action.

Consider the Banach space $\CB(G)$ and let $f\in \CB(G)$ and $\mu \in \M(G)$. Take a net $(g_\alpha)_\alpha$ in $G$ such that $\lim_\alpha g_\alpha = e$ and compute that
    \begin{align*}
        \lim_\alpha (\mu*f)(g_\alpha) & = \lim_\alpha \int_G f(g^{-1}g_\alpha) d\mu(g)              \\
                                      & = \int_G \lim_\alpha  f(g^{-1}g_\alpha) d\mu(g)             \\
                                      & = \int_G f(g^{-1}) d\mu(g) = (\mu*f)(e).
    \end{align*}    
Therefore, $\mu*f \in \CB(G)$.

For the case $E = \RUC(G)$. Take $f \in \RUC(G)$ and $\mu \in \M(G)$. Then $f = \phi*F$ for $\phi \in L^1(G)$ and $F\in L^\infty(G)$ by Cohen-Hewitt Factorization Theorem. This means that $\mu * f = \mu * \phi * F$. But now $\mu*\phi \in L^1(G)$ as this last space can be viewed as the ideal $\M_a(G)$ of $\M(G)$. Employing once again the Cohen-Hewitt Factorization Theorem, we can conclude that $\mu * f \in \RUC(G)$. 

Let's move on the Banach space $\LUC(G)$. Take $f\in \LUC(G)$ and $\mu\in \M(G)$. Then there are $F\in L^\infty(G)$ and $\phi\in L^1(G)$ such that $f = F*\phi$ by Cohen-Hewitt Factorization Theorem. Therefore, 
    \begin{align*}
        \mu*f = \mu*(F*\phi) = (\mu*F)*\phi \in \LUC(G)
    \end{align*}
as $\mu*F \in L^\infty(G)$.

For the space $E = \CU(G)$, we can use the Cohen-Hewitt Factorization Theorem in the same way we did for the case $\RUC(G)$.
\end{proof}

The \textbf{support} $\supp(\mu)$ of a measure $\mu \in \M (G)$ is defined as the complement of the largest open subset of $\mu$-measure zero. It follows that $\supp(\mu)$ is the smallest closed set whose complement has measure zero under $\mu$ and that $\supp(\mu) = \supp(|\mu|)$.

Write $\M_{00}(G)$ for the set of all \textbf{signed finite regular measures on $G$ with compact support}. Then $\M_{00}(G)$ is a normed subalgebra of the Banach algebra $\M (G)$. Using the map $f\mapsto \mu_f$ as in Notation \ref{notation map m_f}, the space $\mathcal{C}_{00}(G)$ equipped with the $||\cdot||_1$- norm isometrically embeds into the normed subalgebra $\M_{00}(G)_a = \M_{00}(G)\cap \M(G)_a$ of all compactly supported signed regular Borel measures which are absolutely continuous with respect to $m_G$. Note that $\M_{00}(G)_a$ is not an ideal in $\M(G)$. However, $\M_{00}(G)_a$ is an ideal in $\M_{00}(G)$. 

\begin{lem}\label{lemma C_00 algebraic ideal}
The algebra $\M_{00}(G)_a$ is an ideal in $\M_{00}(G)$. In particular, the function space  $\mathcal{C}_{00}(G)$ is $\M_{00}(G)$-invariant.
\end{lem}
\begin{proof}
We start showing that $\M_{00}(G)_a$ is an ideal in $\M_{00}(G)$. Let $\mu_1 \in \M_{00}(G)$ and $\mu_2 \in \M_{00}(G)_a$. Then $\mu_1*\mu_2 \in \M(G)_a$. There is only to show that $\mu*f$ has compact support. But this is straightforward as
    \begin{align*}
        \supp(\mu_1*\mu_2) \subset \supp(\mu_1)\supp(\mu_2)
    \end{align*}
by \cite[Chap.VIII §1 No.4 Proposition 5 a)]{bourbakiintegral}. Now the set $\supp(\mu_1)\supp(\mu_2)$ is compact because the multiplication of two compact subsets of a topological group is still a compact subset (\cite[Theorem (4.4)]{abstractharmonic}). This implies that $\supp(\mu_1* \mu_2 )$ is compact as it is closed. Therefore, $\mu_1*\mu_2 \in \M_{00}(G)_a$. In the same way, it is possible to show that $\mu_2*\mu_1 \in \M_{00}(G)_a$. We can conclude that $\M_{00}(G)_a$ is an ideal in $\M_{00}(G)$.

Let now $\mu\in \M_{00}(G)$ and $\phi \in \mathcal{C}_{00}(G)$. We know that $\mathcal{C}_{00}(G) \subset \CU(G)$. Therefore, $\mu*\phi$ is in $\CU(G)$. As before, there is only to show that $\mu*f$ has compact support. But
    \begin{align*}
        \supp(\mu*f) = \supp(\mu*\mu_f) \subset \supp(\mu)\supp(\mu_f) = \supp(\mu)\supp(f).
    \end{align*}
Therefore, $\mathcal{C}_{00}(G)$ is $\M_{00}(G)$-invariant.
\end{proof}

\subsection{Domination and measures} The goal now is to apply measure theory to generalize the theory developed in Section \ref{Chapter Ddominating and majorizing} to locally compact groups. Note that when $G$ is a discrete group, then everything done in this subsection coincides with what was defined in Section \ref{Chapter Ddominating and majorizing}. In fact, we want to develop specific tools to use in the case of non-discrete locally compact groups.

\begin{defn}
For a non-zero function $f\in E_+$, we define
    \begin{align*}
        (E,f)_{\M} &:= \Big\{ h \in E : \exists \lambda \in \M_{00}(G)_+ \text{ s.t. } |h| \leq \lambda * f \Big\}.\\
    \end{align*}
\end{defn}

In other words, the set $(E,f)_{\M}$ is the space of all functions of $E$ which are $\M_{00}(G)_+$-dominated by $f$.

\medskip
By definition, the space $(E,f)_{\M}$ is an ideal of $E$ for every non-zero $f\in E_+$. In particular, $(E,f)_{\M}$ is a Riesz subspace of $E$.

\medskip
If $G$ is a discrete group, then $\M_{00}(G)_+ = c_{00}(G)_+$. Therefore, $(E,f)_\M = (E,f)$ for every non-zero $f\in E_+$. If $G$ is a non-discrete locally compact group, we could not show, either disproved, this last equality. However, the relation between $(E,f)$ and $(E,f)_\M$ is described in the following lemma.

\begin{lem}\label{lemma two compactly supported functions mutually g-domianted}
Let $f\in E_+$ be a non-zero function. Then
    \begin{itemize}
        \item[a)] $(E, f) \subset (E,f)_\M$ and $(E, \phi*f) \subset (E,f)_\M\, $ for every non-zero $\phi\in \mathcal{C}_{00}(G)_+$;
        
        \item[b)] $(E,\phi*f) = (E,\tilde{\phi}*f)$ for every non-zero $\phi,\tilde{\phi}\in \mathcal{C}_{00}(G)_+$;
        
        \item[c)] $(E,\phi * f) = (E,\phi*f)_\M $ for every non-zero $\phi\in \mathcal{C}_{00}(G)_+$;
        
        \item[d)] $(E,\phi*f)_\M = (E,\tilde{\phi}*f)_\M \subset (E,f)_\M$ for every non-zero $\phi, \tilde{\phi}\in \mathcal{C}_{00}(G)_+.$
    \end{itemize}
\end{lem}
\begin{proof}
We start by proving point a). As every positive finite combination of Dirac masses is in $\M_{00}(G)_+$, the inclusion $(E,f)\subset (E,f)_{\M}$ holds. Similarly for every non-zero $\phi\in \mathcal{C}_{00}(G)_+$, the inclusion $(E, \phi*f) \subset (E,f)_\M\, $ also holds because $\mathcal{C}_{00}(G)$ can be seen as a subspace of $\M_{00}(G)$.

The proofs of points b) and c) are only a consequence of the fact that two non-zero positive compactly supported continuous functions always $G$-dominate each other, see point 2) of Example \ref{example invariant normalized integral property}. Indeed, fix two non-zero functions $\phi$ and $\tilde{\phi}$ in $\mathcal{C}_{00}(G)_+$. It suffices to show that there are elements $g_1,...,g_n, x_1,...,x_m\in G$ such that $\phi*f \leq \sum_{j=1}^n g_j\left( \tilde{\phi}*f \right)$ and $\tilde{\phi}*f \leq \sum_{k=1}^m x_k\left(\phi*f \right)$. As $\phi$ and $\tilde{\phi}$ are compactly supported continuous functions, there are $g_1,...,g_n\in G$ such that $\phi\leq \sum_{j=1}^ng_j\tilde{\phi}$ and $x_1,...,x_m\in G$ such that $\tilde{\phi}\leq \sum_{k=1}^m x_j\phi$. Therefore,
    \begin{align*}
        \phi*f \leq \left(\, \sum_{j=1}^n g_j\tilde{\phi} \right)*f = \sum_{j=1}^n g_j\left( \tilde{\phi}*f \right)
        \quad\text{ and }\quad \tilde{\phi}*f \leq \left( \sum_{k=1}^m x_k\phi \right)*f = \sum_{k=1}^m x_k\left(\phi*f \right).
    \end{align*}
The proof to show that $(E, \phi*f)_\M =(E,\tilde{\phi}*f)_\M$ is similar. 
    
Let's prove point c). Let $h\in (E,\phi*f)$. Then $h\in (E,\phi*f)_\M$ by point a). Therefore, $(E,\phi*f) \subset (E,\phi*f)_\M$.  Conversely, let $h\in (E,\phi*f)_\M$. Then there is $\mu\in \M_{00}(G)_+$ such that $|h|\leq \mu* (\phi*f) = (\mu*\phi)*f$. But $\mu*\phi \in \mathcal{C}_{00}(G)$ as it can be seen as an (algebraic) ideal in $\M_{00}(G)$ (Lemma \ref{lemma C_00 algebraic ideal}). So, $h\in (E, \mu*\phi*f) = (E, \phi*f)$ by point b). We can conclude that $(E,\phi*f) = (E,\phi*f)_\M$.

Let's look at point d). Let $\phi,\tilde{\phi}\in \mathcal{C}_{00}(G)_+$ be non-zero. Then
    \begin{align*}
        (E,\phi*f)_\M = (E,\phi*f) = (E,\tilde{\phi}*f) = (E, \tilde{\phi}*f)_\M
    \end{align*}
by points b) and c). Finally,
    \begin{align*}
        (E,\phi*f)_\M = (E,\phi*f) \subset (E,f)_\M
    \end{align*}
by point a).
\end{proof}

It is natural to extend the notion of dominated norm also for spaces of the form $(E,f)_\M$.

\begin{defn}
For a non-zero function $f\in E_+$, we define the map $\overp_f$ on $(E,f)_\M$ as 
     \begin{align*}
          \overp_f(h) := \inf \Big\{ ||\lambda||_{\TV} : |h|\leq \lambda *f \text{ for } \lambda \in \M_{00}(G)_+ \Big\}.
     \end{align*}
\end{defn}

Just as in Section \ref{Chapter Ddominating and majorizing}, we have the following proposition:

\begin{prop}\label{propriety of p_M}
Let $f\in E_+$ be a non-zero function. Then:
    \begin{itemize}
        \item[a)] the pair $\left((E,f)_\M, \overline{p}_f \right)$ is a normed Riesz space on which $G$ acts by positive linear isometries;
        \item[b)] the inequality $||h||_\TV \leq \overp_f(h)||d||_\TV$ holds for every $h\in (E,f)_\M$;
        \item[c)] let $h,v\in (E,f)_\M$ such that $h$ is $\M_{00}(G)_+$-dominated by $v$ and $v\geq 0$, then
                    \begin{align*}
                        \overp_f(h)\leq \overp_v(h)\overp_f(v).
                    \end{align*}
    \end{itemize}
\end{prop}

The proof of this last proposition is quite identical to the proof given for Proposition \ref{prop proprieta $G$-dominated norm}. Therefore, we decided to skip it. 

\medskip
We continue by comparing the $\overp_f$-norm and the $p_f$-norm.

\begin{prop}\label{Convergenza in p_f implica in p_M^f}
Let $f\in E_+$ be a non-zero function and let $(h_\alpha)_\alpha \subset (E,f)_\M$ be a net which converges to $h\in (E,f)_\M$ in $p_F$-norm, where $F$ is  a non-zero positive element of $(E,f)_\M$. Then $(h_\alpha)_\alpha$ also converges to $h$ in $\overp_f$-norm.
\end{prop}

\begin{proof}
As the net $(h_\alpha)_\alpha$ converges to $h$ in $p_F$-norm, for every $\alpha$ there is $n_\alpha \in \N$ and elements $g_1^\alpha,...,g_{n_\alpha}^\alpha\in G$ and $t_1^\alpha, ...,t^\alpha_{n_\alpha} \in \R_+$ such that
    \begin{align*}
        |h_\alpha- h| \leq \sum_{j=1}^{n_\alpha}t_j^\alpha g_j^{\alpha}F \quad \text{ and }\quad  \lim_\alpha \sum_{j=1}^{n_\alpha}t_j^\alpha = 0.
    \end{align*}
Take $\lambda  \in \M_{00}(G)_+$ such that $|F| \leq \lambda*f$. Then
    \begin{align*}
         |h_\alpha- h| \leq \sum_{j=1}^{n_\alpha}t_j^\alpha g_j^{\alpha}(\lambda*f) = \underbrace{\left(\,\sum_{j=1}^{n_\alpha}t_j^\alpha \delta_{g_j^{\alpha}}\right)*\lambda}_{\in \M_{00}(G)_+} *f \quad \text{ for every }\alpha. 
    \end{align*}
Therefore,
    \begin{align*}
        \lim_\alpha \Big|\Big| \left(\,\sum_{j=1}^{n_\alpha}t_j^\alpha \delta_{g_j^{\alpha}}\right)*\lambda \Big|\Big|_\TV \leq \lambda(G) \lim_\alpha \left(\, \sum_{j=1}^{n_\alpha}t_j^\alpha \delta_{g_j^{\alpha}}\right)(G)  \leq \lambda(G) \lim_\alpha  \sum_{j=1}^{n_\alpha}t_j^\alpha = 0,
    \end{align*}
which proves that $(h_\alpha)_\alpha$ convergences to $h$  in $\overp_f$-norm.
\end{proof}

In particular, the $\overline{p}_f$-norm is weaker than the $p_f$-norm on the space $(E,f)$. A situation where the two norms are actually equal is given in the following proposition. 

\begin{prop}
Let $f\in E_+$ be a non-zero function. Then 
    \begin{align*}
        \left( (E, \phi*f), p_{\phi*f} \right) = \left( (E,\phi*f)_\M, \overp_{\phi*f} \right)
    \end{align*}
for every non-zero $\phi\in \mathcal{C}_{00}(G)_+$.
\end{prop}
\begin{proof}
By Lemma \ref{lemma two compactly supported functions mutually g-domianted}, we can fix an arbitrary non-zero $\phi \in \mathcal{C}_{00}(G)_+$. We know that $(E,\phi*f) = (E,\phi*f)_\M$ and that $\overp_{\phi*f} \leq p_{\phi*f}$ by  point c) of Lemma \ref{lemma two compactly supported functions mutually g-domianted} and Proposition \ref{Convergenza in p_f implica in p_M^f}, respectively. Therefore, it suffices to show that $p_{\phi*f} \leq \overp_{\phi*f}$. To this aim, let $h\in (E, \phi*f)_\M$. Then there is a net $(\mu_\alpha)_\alpha $ in $\M_{00}(G)_+$ such that $|h|\leq \mu_\alpha*\phi*f$ for every $\alpha$ and $\lim_\alpha ||\mu_\alpha||_\TV = \overp_{\phi*f}(h)$. Thus,
    \begin{align*}
        p_{\phi*f}(h) & \leq p_{\phi*f}(\mu_\alpha*\phi*f)
                    = p_{\phi*f}\left( \int_G x(\phi*f)d\mu_\alpha(x) \right) \\
                    & \leq \int_G p_{\phi*f}(x(\phi*f))d\mu_\alpha(x)
                    = ||\mu_\alpha||_\TV
    \end{align*}
for every $\alpha$. Hence, $p_{\phi*f}(h) \leq \overp_{\phi*f}(h)$.
\end{proof}

For dominating norms of the form $p_{\phi*f}$ there is also an interesting convergence property illustrated in the following lemma. This result will be fundamental in proof of Theorem \ref{Equivalence of integral on different spaces}.

\begin{lem}\label{lemma app identity}
Let $f\in E_+$ and $\phi\in \mathcal{C}_{00}(G)_+$ be non-zero functions. Then for every $\theta\in \mathcal{C}_{00}^1(G)$ and for every bounded approximate identity $(e_\alpha)_\alpha$ for $L^1(G)$ in $\mathcal{C}_{00}^1(G)$ with decreasing support, we have that 
    \begin{align*}
        \lim_\alpha \theta*e_\alpha*h = \theta*h
        \quad \text{ and } \quad \lim_\alpha e_\alpha*\theta*h = \theta*h
    \end{align*}
in $p_{\phi*f}$-norm for every $h\in (E,f)_\M$ and every $\theta\in \mathcal{C}_{00}(G)$.
\end{lem}

The existence of bounded approximate identities with decreasing support is assured by Urysohn Lemma and the fact that $G$ is locally compact.

\begin{proof}[Proof of Lemma \ref{lemma app identity}]
Let $h\in (E,f)_\M$, $\theta\in \mathcal{C}_{00}(G)$ and $(e_\alpha)_\alpha$ a bounded approximate identity for $L^1(G)$ in $\mathcal{C}_{00}^1(G)$ with decreasing support. We want to show that $\lim_\alpha \theta*e_\alpha*h = \theta*h$ in $p_{\phi*f}$-norm. Without loss of generality, we can suppose that $h$ is positive, since $(E,f)_\M$ is spanned by positive elements. Now the estimation
    \begin{align*}
        p_{\phi*f}(\theta*e_\alpha*h - \theta*h )&  = 
            p_{\phi*f}\big((\theta*e_\alpha - \theta)*h \big)
            \leq p_{\phi*h}\big((\theta*e_\alpha - \theta)*h\big)p_{\phi*f}(\phi*h)
    \end{align*}
holds for every $\alpha$. The last inequality is possible thanks to point c) of Proposition \ref{prop proprieta $G$-dominated norm}. The second term is finite because $\phi*h \in(E,\phi*f)_\M$. So let's focus on the first term. By Proposition \ref{prop conv of net in sup norm implies convergence in p_f}, $\lim_\alpha \theta*e_\alpha = \theta$ in $p_\phi$-norm. Therefore, for every $\alpha$ there is $n_\alpha\in \N$ and elements $g_1^{\alpha},...,g_{n_\alpha}^\alpha \in G$ and $t_1^\alpha,...,t_{n_\alpha}^\alpha\in \R_+$ such that
    \begin{align*}
        |\theta*e_\alpha- \theta| \leq \sum_{j=1}^{n_\alpha}t_j^\alpha g_j^{\alpha}\phi \quad \text{ and }\quad  \lim_\alpha \sum_{j=1}^{n_\alpha}t_j^\alpha = 0.
    \end{align*}
This implies that
    \begin{align*}
        |\theta*e_\alpha- \theta|*h \leq 
        \left(\,\sum_{j=1}^{n_\alpha}t_j^\alpha g_j^{\alpha}\phi\right)*h
        = \sum_{j=1}^{n_\alpha}t_j^\alpha g_j^{\alpha}(\phi*h).
    \end{align*}
Hence,
    \begin{align*}
        p_{\phi*h} \left( |\theta*e_\alpha- \theta|*h \right)  \leq p_{\phi*h} \left( \; \sum_{j=1}^{n_\alpha}t_j^\alpha g_j^{\alpha}(\phi*h) \right) 
            \leq \sum_{j=1}^{n_\alpha}t_j^\alpha p_{\phi*h}(g_j^{\alpha}(\phi*h))
            =  \sum_{j=1}^{n_\alpha}t_j^\alpha,
    \end{align*}
for every $\alpha$. Taking the limit with respect to $\alpha$ of this last inequality, it follows that $\lim_\alpha \theta*e_\alpha*h = \theta*h$ in $p_{\phi*f}$-norm. 

The proof to show that $\lim_\alpha e_\alpha*\theta*h = \theta*h$ in $p_{\phi*f}$-norm is similar. 
\end{proof}

\medskip
We don't know if the natural action of $G$ on the normed Riesz space $\left(\RUC(G,f)_\M, \overp_f \right)$ is orbitally continuous. However, it is possible to prove something sufficient for our purposes.

\begin{lem}\label{Orbitally continuoity for subspace D}
Let $f\in E_+$ be a non-zero function and define the $\M_{00}(G)$-invariant linear subspace
    \begin{align*}
        D := \Span_\R \Big\{ \psi*h : \psi \in \mathcal{C}_{00}(G) \text{ and } h \in (E,f)_\M \Big\} \subset (E,f)_\M.
    \end{align*}
Then the following assertions hold: 
\begin{itemize}
    \item[a)] the action of $G$ on $D$ is orbitally continuous with respect to the $p_{\phi*f}$-norm for every non-zero $\phi \in \mathcal{C}_{00}(G)_+$;
    \item[b)] the action of $G$ on $D$ is orbitally continuous with respect to the $\overp_f$-norm.
\end{itemize}
\end{lem}

\begin{proof}
First of all, we should check that $D$ is really a $\M_{00}(G)$-invariant linear subspace of $(E,f)_\M$. By definition, $D$ naturally carries a structure of linear space. Therefore, we only need to show that it is a subspace of $(E,f)_\M$. To this aim, let $v\in D$. Then there are $\psi\in \mathcal{C}_{00}(G)$ and $h\in(E,f)_\M$ such that $v = \psi*h$. As $h\in(E,f)_\M$, there is $\tau \in \M_{00}(G)_+$ such that $|h|\leq \tau*f$. Thus,
    \begin{align*}
        |v|\leq |\psi*h| = |\lambda_\psi * h| \leq |\lambda_\psi|*|h| \leq \underbrace{\big(|\lambda_\psi|*\tau\big)}_{\in  \M_{00}(G)_+}*f
    \end{align*}
showing that $D$ is a subspace of $(E,f)_\M$. The $\M_{00}(G)$-invariance is implied by the definition of $D$ and because the space $\mathcal{C}_{00}(G)$ can be identified with the ideal $\M_{00}(G)_a$.

We proceed proving point a). Fix $ \psi*h \in D$ and let $(g_\alpha)_\alpha \subset G$ be a net such that $\lim_\alpha g_\alpha = e$, the identity element of $G$. Without loss of generality, we can suppose that $\psi$ is positive, since each element in $\mathcal{C}_{00}(G)$ is the difference of two positive elements. Then the estimation
    \begin{align*}
        p_{\phi*f}(g_\alpha \big(\psi*h )- \psi*h ) & = p_{\phi*f}( (g_\alpha \psi- \psi)*h ) \\
         & \leq p_{\phi*f}( |g_\alpha \psi- \psi|*|h| ) \\  & \leq p_{\psi*|h|}(|g_\alpha \psi- \psi|*|h|)p_{\phi*f}(\psi*|h|)
    \end{align*}
holds for every $\alpha$. The last inequality is possible thanks to point c) of Proposition \ref{prop proprieta $G$-dominated norm}. Now $p_{\phi*f}(\psi*|h|) < \infty$ because $\psi*|h| \in (E,\phi* f)$. Thus, let's look at the first term. To this aim, note that the action of $G$ on $\mathcal{C}_{00}(G)$ is orbitally continuous with respect to the $p_\psi$-norm thanks to Lemma \ref{la azione é orbitalli continuous}. This means that, for every $\alpha$, there is $n_\alpha \in \N$ and elements $t_1^{\alpha},...,t_{n_\alpha}^{\alpha}\in \R_+$ and $g_1^{\alpha},...,g_{n_\alpha}^\alpha \in G$ such that
    \begin{align*}
        |g_\alpha\psi - \psi | \leq \sum_{j=1}^{n_\alpha}t_{j}^\alpha g_j^\alpha \psi \quad \text{ and } \quad \lim_{\alpha}\sum_{j=1}^{n_\alpha}t_j^{\alpha} = 0.
    \end{align*}
Taking the convolution with $|h|$ on both sides of the previous inequality, we have that 
    \begin{align*}
        |g_\alpha\psi - \psi |*|h| \leq \left(\, \sum_{j=1}^{n_\alpha}t_{j}^\alpha g_j^\alpha \psi\right)*|h| = \sum_{j=1}^{n_\alpha}t_{j}^\alpha g_j^\alpha \left( \psi*|h|\right).
    \end{align*}
This led to the estimation
    \begin{align*}
        p_{\psi*|h|}(|g_\alpha\psi - \psi |*|h|)& \leq 
             p_{\psi*|h|}\left(\, \sum_{j=1}^{n_\alpha}t_{j}^\alpha g_j^\alpha \left( \psi*|h|\right) \right)
             \leq \sum_{j=1}^{n_\alpha}t_{j}^\alpha  p_{\psi*|h|}\left(g_j^\alpha ( \psi*|h|)\right)
             = \sum_{j=1}^{n_\alpha}t_{j}^\alpha,
    \end{align*}
which implies that $\lim_\alpha p_{\psi*|h|}(|g_\alpha \psi- \psi|*|h|) = 0 $. We can hence conclude that the action is orbitally continuous.

The proof of point b) is straightforward using Proposition \ref{Convergenza in p_f implica in p_M^f} and point a).
\end{proof}

\subsection{A new kind of integral} The adaption of dominating spaces to the measurable setting gives rise to a new notion of normalized integral and consequently a new normalized integral property. 

Define the two sets
\begin{align*}
        \M^1_{00}(G) := \left\{ \lambda \in \M_{00}(G): \lambda \geq 0 \text{ and }||\lambda||_\TV = \lambda(G)=1 \right\}
\end{align*}
and 
\begin{align*}
        \mathcal{C}^1_{00}(G) := \left\{ \phi \in \mathcal{C}_{00}(G): \phi \geq 0 \text{ and }||\phi||_1 = 1 \right\}.
\end{align*}

Then $\M^1_{00}(G)$ and $\mathcal{C}_{00}^1(G)$ are semigroups under convolution.\footnote{We mean convolution between measures for the first one, while for the second one, convolution between functions.}

\medskip
Let $f\in E_+$ be a non-zero vector. A functional $\psi \in (E,f)_\M^*$ is said \textbf{measurably invariant}, or \textbf{$\M^1_{00}(G)$-invariant}, if $\psi( \lambda * v) = \psi(v)$ for every $\lambda \in \M^1_{00}(G)$ and every $v\in (E,f)_\M$.  

\begin{defn}
 We say that a locally compact group $G$ has \textbf{the measurably invariant normalized integral property for $E$} if there exists a positive linear $\M_{00}^1(G)$-invariant functional $\I$ on $(E,f)_\M$ such that $\I(f)=1$ for every non-zero  $f\in E_+$.
\end{defn}

Note that the notions of measurably invariant normalized integral and invariant normalized integral coincide for discrete groups. However, the measurably invariant normalized integral property is a priori stronger than the invariant normalized integral property for non-discrete locally compact groups. In fact, 

\begin{prop}\label{measured invariant on E implies invariant on E}
Suppose that $G$ has the measurably invariant normalized integral property for $E$, then $G$ has the invariant normalized integral property for $E$.
\end{prop}

\begin{proof}
The proof is directly deduced from the fact that $(E,f) \subset (E,f)_\M$  for every non-zero positive $f\in E$, and the fact that $G$ can be seen as a subgroup of $\M_{00}(G)$ via the map $g \longmapsto \delta_g$.
\end{proof}

Finally, measurably invariant normalized integrals are $\overp_f$-continuous. 

\begin{prop}\label{continuoity of I}
Let $f\in E_+$ be a non-zero function, and let $\I$ be a measurably invariant normalized integral defined on $(E,f)_\M$. Then $\I$ is continuous with respect to the $\overp_f$-norm and has norm operator equal to one.
\end{prop}

\begin{proof}
On the one hand, $\I(f)=1$ and $\overp_f(f)=1$. Thus, $||\I||_{op}\geq 1$. On the other hand, let $h\in (E,f)_\M$ be such that $\overp_f(h) = 1$. Then there is a net $(\lambda_\alpha)_\alpha \subset \M_{00}(G)_+$ such that $|h|\leq \lambda_\alpha*f$ for every $\alpha$ and such that $\lim_\alpha ||\lambda_\alpha||_\TV = 1$. Therefore, 
    \begin{align*}
        |\I(h)| \leq \I(\lambda_\alpha*f) = ||\lambda_\alpha ||_\TV \; \I \left(\frac{\lambda_\alpha}{||\lambda_\alpha||_\TV}*f \right) = ||\lambda_\alpha||_\TV \quad \text{ for every } \alpha.
    \end{align*}
Hence, $|\I(h)|\leq 1$ which implies $||\I||_{op}\leq 1$.
\end{proof}

\subsection{Proof of Theorem \ref{Equivalence of integral on different spaces}} The entire subsection is dedicated to the demonstration of Theorem \ref{Equivalence of integral on different spaces}. We start by stating and proving a small lemma that solve the technical core of the theorem. After that, two propositions prove intermediate results. The proof of the theorem can be founded at the end of the subsection. 

\begin{lem}\label{lemma restriction of I measurably invariant}
Let $f\in \RUC(G)_+$ be a non-zero function and let $\phi \in \mathcal{C}_{00}^1(G)$. Consider the subspace
    \begin{align*}
        D = \Span_\R \left\{ \psi * h : \psi \in \mathcal{C}_{00}(G) \text{ and } h \in \RUC(G, f)_\M \right\} \subset \RUC(G,\phi*f).
    \end{align*}
Then every positive $G$-invariant $p_{\phi*f}$-continuous functional $\I$  on $D$ is also measurably invariant, i.e., $\I(\mu * h) = \I(h)$ for every $\mu\in \M_{00}^1(G)$ and $h\in D$. 
\end{lem}

\begin{proof}
First of all, note that the action of $G$ on the closure $\overline{D}^{p_{\phi*f}}\subset \RUC(G,\phi*f)$ is orbitally continuous with respect to the $p_{\phi*f}$-norm because of point a) of Lemma \ref{Orbitally continuoity for subspace D} and \cite[Lemma $4.1.9$]{max}. Moreover, the functional $\I$ extends uniquely to a positive $G$-invariant linear functional on $\overline{D}^{p_{\phi*f}}$, since it is positive, $G$-invariant and $p_{\phi*f}$-continuous on $D$. 
Now take $h\in \overline{D}^{p_{\phi*f}}$ and $\mu\in \M_{00}^1(G)$ and consider the function
    \begin{align*}
        F : (G, \mu) \longrightarrow \overline{D}^{p_{\phi*f}}, \quad g \longmapsto F(g) = gh.
    \end{align*}
The map $F$ is Bochner integrable, because it is continuous and the real integral
    \begin{align*}
        \int_G p_{\phi*f}(F(g))d\mu(g) = \int_G p_{\phi*f}\big(gh\big) d\mu(g) \leq p_{\phi*f}(h) \mu(G)
    \end{align*}
is finite. Therefore, the Bochner integral $\int_G F(g) \mu(g) \in \overline{D}^{p_{\phi*f}}$, see \cite[Appendix B]{principle}.
The $p_{\phi*f}$-continuity of $\I$ on $\overline{D}^{p_{\phi*f}}$ gives the following:  
    \begin{align*}
        \I(\mu*h)  & = \I\left(\int_G gh \,d\mu(g) \right)       
                             = \int_G \I\left( gh \right) d\mu(g)      \\
                            & = \int_G \I\left(h\right) d\mu(g)           
                             = \I(h)\mu(G) = \I(h).
    \end{align*}
This proves that $\I$ is also measurably invariant on $D$.

\end{proof}

We now prove the above-cited two propositions.

\begin{prop}\label{measured invariant on L implies measured invariant on RUC}
Suppose that $G$ has the measurably invariant normalized integral property for $\RUC(G)$, then $G$ has the measurably invariant normalized integral property for $L^\infty(G)$.
\end{prop}

\begin{proof}
Let $f\in L^\infty(G)$ be a non-zero positive function and let $\phi \in \mathcal{C}_{00}^1(G)$. Define the linear operator
    \begin{align*}
        T: L^\infty(G,f)_\M \longrightarrow \mathcal{C}_{ru}^b(G,\phi*f)_\M,\quad h\longmapsto \phi*h.
    \end{align*}
First of all, we check that $T$ is well-defined. It is clear that $\phi*h\in\mathcal{C}_{ru}^b(G)$ for every $h\in L^\infty(G)$ by Cohen-Hewitt Factorization Theorem. Suppose that $h\in L^\infty(G,f)_\M$, then there is $\lambda \in \M_{00}(G)_+$ such that $|h|\leq \lambda*f$. Consequently,
    \begin{align*}
        |\phi*h| \leq \phi*|h| \leq \underbrace{\phi*\lambda}_{\in\, \mathcal{C}_{00}(G)}*f,
    \end{align*}
and so $\phi*h \in \mathcal{C}_{ru}^b(G,(\phi*\lambda)*f)_\M = \RUC(G,\phi*f)_\M$ by point d) of Lemma \ref{lemma two compactly supported functions mutually g-domianted}.

Let now $\I$ be a measurably invariant normalized integral on $\mathcal{C}_{ru}^b(G,\phi*f)_\M$ and define on $L^\infty(G,f)_\M$ the functional $\overline{\I} = \I \circ T$. We claim that $\overline{\I}$ is a measurably invariant normalized integral for $L^\infty(G,f)_\M$. Clearly, $\overline{\I}$ is linear as a composition of linear maps and it is normalized because $\overline{\I}(f)= \I(\phi*f)=1$. Therefore, it is left to show that $\overline{\I}$ is $\M_{00}^1(G)$-invariant. 
To this end, let $\mu\in \M_{00}^1(G)$, $h\in L^\infty(G,f)_\M$ and take a bounded approximate identity $(e_\alpha)_\alpha$ for $L^1(G)$ in $\mathcal{C}_{00}^1(G)$ with decreasing support. Then
    \begin{align*}
        \overline{\I}(\mu*h) & = \I ( \phi* \mu *h)
            = \I\left( \lim_\alpha (\phi*\mu)*e_\alpha *h\right)
            = \lim_\alpha \I\left( (\phi*\mu)*e_\alpha*h \right) \\
            & = \lim_\alpha \I\left( \phi*e_\alpha*h\right)
            = \I\left( \lim_\alpha \phi*e_\alpha*h\right)
            = \I \left( \phi*h \right) = \overline{\I}(h),
    \end{align*}
since $\lim_\alpha (\phi*\mu)*e_\alpha *h = (\phi*\mu)*h$ and $\lim_\alpha \phi*e_\alpha*h = \phi*h$ in $p_{\phi*f}$-norm by Lemma \ref{lemma app identity}.
This shows that $\overline{\I}$ is measurably invariant and concludes the proof.
\end{proof}

\begin{prop}\label{invariant for RUC implies measured invariant fo RUC}
Suppose that $G$ has the invariant normalized integral property for $\mathcal{C}_{ru}^b(G)$, then $G$ has the measurably invariant normalized integral property for $\mathcal{C}_{ru}^b(G)$.
\end{prop}

\begin{proof}
Let $f\in \mathcal{C}_{ru}^b(G)_+$ be a non-zero function and chose a $\phi\in \mathcal{C}_{00}^1(G)$.

Consider the linear operator
    \begin{align*}
        T: \RUC(G,f)_\M \longrightarrow \RUC(G,\phi*f), \text{ }h \longmapsto \phi*h
    \end{align*}
and note that it is well-defined because if $h\in \mathcal{C}_{ru}^b(G,f)_\M$, then there is $\mu\in \M_{00}(G)_+$ such that $|h|\leq \mu*f$ which implies that
    \begin{align*}
        |\phi*h| \leq \phi*|h| \leq \underbrace{(\phi*\mu)}_{\in\, \mathcal{C}_{00}(G)}*f.
    \end{align*}
Hence, $\phi*h \in \RUC(G, (\phi*\mu)*f)=\RUC(G, \phi*f)$ by point b) of Lemma \ref{lemma two compactly supported functions mutually g-domianted}.

Moreover, note that $T$ maps $\RUC(G,f)_\M$ into 
    \begin{align*}
        D = \Span_\R \left\{ \psi * h : \psi \in \mathcal{C}_{00}(G) \text{ and } h \in \RUC(G, f)_\M \right\} \subset\RUC(G,\phi*f).
    \end{align*}

By hypothesis, there is an invariant normalized integral $\I$ on $\RUC(G,\phi*f)$. We define the linear functional $\overline{\I}$ on $\RUC(G,f)_\M $ as $\overline{\I} = \I \circ T$, and we claim that it is a measurably invariant normalized integral. Indeed, $\overline{\I}$ is normalized because $\overline{\I}(f) = \I(\phi*f)=1$ and linear as composition of linear maps. Thus, it is left to show that $\overline{\I}$ is measurably invariant. To this end, let $\mu\in \M_{00}^1(G)$, $h\in \RUC(G,f)_\M$ and take a bounded approximate identity $(e_\alpha)_\alpha$ for $L^1(G)$ in $\mathcal{C}_{00}^1(G)$ with decreasing support. Note that $\I$ is measurably invariant when restricted to $D$ by Lemma \ref{lemma restriction of I measurably invariant}. Then
    \begin{align*}
        \overline{\I}(\mu*h) & = \I ( \phi* \mu *h)
            = \I\left( \lim_\alpha (\phi*\mu)*e_\alpha *h\right)
            = \lim_\alpha \I\left( (\phi*\mu)*e_\alpha*h \right) \\
            & = \lim_\alpha \I\left( \phi*e_\alpha*h\right)
            = \I\left( \lim_\alpha \phi*e_\alpha*h\right)
            = \I \left( \phi*h \right) = \overline{\I}(h),
    \end{align*}
since $\lim_\alpha (\phi*\mu)*e_\alpha *h = (\phi*\mu)*h$ and $\lim_\alpha \phi*e_\alpha*h = \phi*h$ in $p_{\phi*f}$-norm by Lemma \ref{lemma app identity}.
This shows that $\overline{\I}$ is measurably invariant and concludes the proof.

\end{proof}

We are finally ready to give the proof of Theorem \ref{Equivalence of integral on different spaces}.

\begin{proof}[Proof of Theorem \ref{Equivalence of integral on different spaces}]
Firstly, if $G$ has the invariant normalized integral property for $L^\infty(G)$, then $G$ has the invariant normalized integral property for all the other function spaces. Secondly, the invariant normalized integral property for $\UC(G)$ implies the invariant normalized integral property for $\RUC(G)$. Indeed, let $f\in \RUC(G)_+$ be a non-zero function. Fix a non-zero $\phi\in \mathcal{C}_{00}(G)$, and define the positive linear operator
    \begin{align*}
        T: \RUC(G,f) \longrightarrow \UC(G, f*\phi), 
        \quad h \longmapsto T(h) = h*\phi.
    \end{align*}
The operator $T$ is well-defined thanks to the Hewitt-Cohen Factorization Theorem. By hypothesis, there is an invariant normalized integral $\I$ on $\UC(G, \phi*f)$. Then the functional $\overline{\I} = \I \circ T$ is an invariant normalized integral for $\RUC(G, f)$.

Therefore, it is enough to show that the invariant normalized integral property for $\RUC(G)$ implies the invariant normalized integral property for $L^\infty(G)$. Thus, suppose that $G$ has the invariant normalized integral property for $\RUC(G)$. Then $G$ has the measurably invariant normalized integral property for $\RUC(G)$ by Proposition \ref{invariant for RUC implies measured invariant fo RUC}. Therefore, $G$ also has the measurably invariant normalized integral property for $L^\infty(G)$ by Proposition \ref{measured invariant on L implies measured invariant on RUC}. We can finally conclude that $G$ has the invariant normalized integral property for $L^\infty(G)$ using Proposition \ref{measured invariant on E implies invariant on E}.
\end{proof}

\begin{rem}
Note that if we fix $f = \1_G$, then we have the case of amenability \cite[Theorem $2.2.1$]{greenleaf}. This is because $\1_G$ is a fixed-point for the action by convolution of the semigroup $\M^1_{00}(G)$ on $E$ and because the space $\mathcal{C}_{00}^1(G)$ is $||\cdot||_1$-dense in 
    \begin{align*}
        \text{P}(G) = \left\{ \phi \in L^1(G) : ||\phi||_1 = 1  \text{ and } \phi \geq 0 \right\}.
    \end{align*}
\end{rem}

\subsection{Variations on the translate property} We extend the notion of translate property for all the function spaces considered in this section. After that, we give a complete answer to Greenleaf's question in the locally compact case. Thanks to this result, it will be possible to affirmatively answer Greenleaf's question for all the classical Banach spaces. 

\begin{defn}
  We say that a locally compact group $G$ has \textbf{the translate property for} $E$ if for every non-zero  $f \in E_+$ whenever
    \begin{align*}
        \sum_{j=1}^nt_jg_jf \geq 0 \quad \implies \quad \sum_{j=1}^nt_j \geq 0,
    \end{align*}
 where $t_1,...,t_n\in \R$ and $g_1,...,g_n\in G$.
\end{defn}

There are essentially two questions that arise spontaneously. The first is understanding if the translate property for $\RUC(G)$ implies the invariant normalized integral property for $\RUC(G)$ in the locally compact case, i.e., the Greenleaf's question. The second is understanding if all the translate properties we defined just above are equivalent.

\medskip
We start providing the proof of the solution to the Greenleaf question for the locally compact case. 

\begin{proof}[Proof of Theorem \ref{theorem TP imples measurably integral}.]
Let $f\in \RUC(G)_+$ be a non-zero function and fix a $\phi\in \mathcal{C}_{00}^1(G)$. Consider the functional given by
    \begin{align*}
        \omega: \Span_\R\left\{ g(\phi*f) : g\in G \right\} \longrightarrow \R, \quad \sum_{j=1}^nt_jg_j(\phi*f) \longmapsto \omega\left(\; \sum_{j=1}^nt_jg_j(\phi*f)\right) = \sum_{j=1}^nt_j.
    \end{align*}
Thanks to the translate property, $\omega$ is well-defined. Moreover, $\omega$ is positive, $G$-invariant and $p_{\phi*f}$-continuous. Now we have the inclusion:
    \begin{align*}
        \Span_\R\left\{ g(\phi*f) : g\in G \right\} \subset \Span_\R \left\{ \psi*h: \psi\in \mathcal{C}_{00}(G)\text{ and }h\in \RUC(G,f)_\M \right\} = D.
    \end{align*}
Note that the action of $G$ on $D$ is orbitally continuous with respect to the $p_{\phi*f}$-norm by point a) of Lemma \ref{Orbitally continuoity for subspace D}. Because $G$ is amenable (Theorem \ref{translate property implies amenability}), we can apply \cite[Theorem 3.2 (j)]{lau} and extend $\omega$ to a positive linear $G$-invariant $p_{\phi*f}$-continuous functional $\overline{\omega}$ defined on all $D$. Note that the functional $\overline{\omega}$ is also measurably invariant by Lemma \ref{lemma restriction of I measurably invariant}. Consider the positive linear operator
    \begin{align*}
        T: \RUC(G,f) \longrightarrow D, \quad
         h \longmapsto \phi*h
    \end{align*}
and define the functional $\I = \overline{\omega}\circ T$ on $\RUC(G,f)$. We claim that $\I$ is an invariant normalized integral. Indeed, $\I$ is positive as composition of positive operators and
    \begin{align*}
        \I(f)= \overline{\omega}(\phi*f) = \omega(\phi*f) = 1.
    \end{align*}
Therefore, it is left to show that $\I$ is invariant. Let $g\in G$, $h\in \RUC(G,f)$ and take a bounded approximate identity $(e_\alpha)_\alpha$ for $L^1(G)$ in $\mathcal{C}_{00}^1(G)$ with decreasing support. Then
    \begin{align*}
        \I(gh) &= \overline{\omega}\left(\phi*(gh)\right)
               = \overline{\omega}\left(\phi_g*h \right)
               = \overline{\omega}\left(\lim_\alpha \phi_g*e_\alpha*h \right) \\
               & = \lim_\alpha \overline{\omega}\left( \phi*e_\alpha*h \right)
               = \overline{\omega}\left( \lim_\alpha \phi*e_\alpha*h \right) 
               = \overline{\omega}(\phi*h) = \I(h),
    \end{align*}
since $\lim_\alpha \phi_g*e_\alpha*h = \phi_g*h$ and $\lim_\alpha \phi*e_\alpha*h = \phi*h$ in $p_{\phi*f}$-norm by Lemma \ref{lemma app identity}. Therefore, $\I$ is invariant, and $G$ has the invariant normalized integral property for $\RUC(G)$.

\end{proof}

\begin{rem}\label{remark TP for function F of the form}
This last theorem shows in particular that only the functions $F\in \RUC(G)$ of the form $\phi*f$, where $\phi \in \mathcal{C}_{00}^1(G)$ and $f\in \RUC(G)$, need to have the translate property to ensure that $G$ has the invariant normalized integral property for $\RUC(G)$. 
\end{rem}

Thanks to this last result, it is also possible to solve Greenleaf's question for all classical Banach spaces. 

\begin{proof}[Proof of Theorem \ref{TP on E equivalent to NII on E}]
If $G$ has the invariant normalized property for $E$, then $G$ has the translate property for $E$. So, let's prove the other direction.

The case $E= \RUC(G)$ has already been solved in Theorem \ref{theorem TP imples measurably integral}.

Suppose that $G$ has the translate property for $E\in \left\{ L^\infty(G), \CB(G) \right\}$. Then $G$ has the translate property for $\RUC(G)$. Therefore, $G$ has the invariant normalized integral property for $\RUC(G)$. This implies that $G$ has the invariant integral property for $E$ by Theorem \ref{Equivalence of integral on different spaces}.

Let's move onto the case $E=\LUC(G)$. We claim that $G$ also has the translate property for $L^\infty(G)$. Indeed, let $f\in L^\infty (G)$ be a non-zero function and let $t_1,...,t_n\in \R$ and $g_1,...,g_n\in G$ such that $\sum_{j=1}^n t_jg_jf \geq 0$. Chose a non-zero positive $\phi\in \mathcal{C}_{00}(G)$. Since taking the convolution with a positive function is a positive operation, the inequality 
\begin{align*}
    \left(\; \sum_{j=1}^n t_jg_jf \right) * \phi  = \sum_{j=1}^n t_jg_j \underbrace{\left( f*\phi \right)}_{\in\; \LUC(G)} \geq 0
\end{align*}
holds. Therefore, $\sum_{j=1}^nt_j\geq 0$. We can hence conclude as in the previous case.

Finally, let $E = \UC(G)$. Then $G$ also has the translate property for $\RUC(G)$ as $f*\phi\in \UC(G)$ for every $\phi\in \mathcal{C}_{00}(G)$ and $f\in \RUC(G)$. Thus, we can conclude as above.
\end{proof}

A combination of Theorems \ref{Equivalence of integral on different spaces} and \ref{TP on E equivalent to NII on E} gives a proof of Theorem \ref{Equivalence of TP on different spaces} which answers the second question we asked ourselves.

\begin{scholium}\label{schoulium 1}
In \cite[(6.42)]{paterson} Paterson defined a \textit{generalized translate property} claiming that it could be a good definition for a topological version of the translate property. One could show that the notion defined by Paterson is equivalent to our translate property. The \textit{hard} part of this equivalence can be shown directly using the notion of measurably invariant integral. 

Another topological translate property was defined and studied by Jenkins under the name of \textit{property $(P)$} in \cite{jenkisfolner}. Moreover, in a later publication, he claimed that this property implies the existence of invariant normalized integrals defined on particular linear subspaces of $L^\infty(G)$, see \cite[Proposition $5$]{jenkinsfixedpoint}. Sadly, the proof is incorrect because it erroneously claimed that $G$ acts continuously on $L^\infty(G)$, which is not generally true. Note that it is possible to show that Jenkins's translate property is equivalent to the translate property and Paterson's property.
\end{scholium}


\section{Equivalent Fixed-point Properties}\label{section equivalent fixed-point properties}

This section aims to get rid of some technical details when working with the fixed-point property for cones for locally compact groups. The results here developed will be helpful later to study the class of locally compact groups which enjoy the fixed-point property for cones.

\subsection{The continuous fixed-point property for cones}
The following lemma shows that it is possible to suppose just continuity in the definition of the fixed-point property for cones for locally compact groups. In particular, the generalization of the fixed-point property for cones proposed by Monod in \cite[Example 38]{monod} and our generalization coincide for locally compact groups. In fact, point b) of the following lemma is exactly the generalization suggested by Monod.

\medskip
A representation of a group on a locally convex vector space is said bounded if every orbit is bounded. 

\begin{lem}\label{lemma continuous fixed-point property for cones}
Let $G$ be a locally compact group. The following assertions are equivalent:
    \begin{itemize}
        \item[ a)] the group $G$ has the fixed-point property for cones;
        \item[b)] every bounded orbitally continuous representation of $G$ on a non-empty weakly complete proper convex cone $C$ in a locally convex vector space $E$ which is of cobounded type has a non-zero fixed-point.
    \end{itemize}
\end{lem}
\begin{proof}
We start showing that a) implies b). Suppose that $G$ has a bounded orbitally continuous representation on a non-empty weakly complete proper cone $C$ in a locally convex vector space $E$ which is of cobounded type. Fix a non-zero vector $v \in C$, and define the positive linear operator
    \begin{align*}
        T: E' \longrightarrow \CB(G), \quad \lambda \longmapsto T(\lambda)
    \end{align*}
where $T(\lambda)(g) = \lambda(gv)$ for every $g\in G$. The operator $T$ is well-defined as composition of continuous maps and $T(\lambda)$ is bounded for every $\lambda\in E'$ because the set $Gv \subset E$ is bounded and continuous linear functionals map bounded sets to bounded sets (\cite[III §1 No.3 Corollaire 1]{bourbakiespace}). Moreover, $T$ is also equivariant with respect to the left-translation representation of $G$ on $\CB(G)$. Now, let $d\in E'$ be the $G$-dominating element given by the cobounded condition. Therefore, $T$ maps $E'$ into $\CB(G, T(d))$. On this last space there is an invariant normalized integral. This is because, for locally compact groups, the fixed-point property for cones of $G$ is equivalent to the invariant normalized integral property for $\CB(G)$ (Theorems \ref{FPC equivalente integral} and \ref{Equivalence of integral on different spaces}). We can now employ the same strategy of the implication $(4) \implies (1)$ of \cite[Theorem 7, pp. 77-78]{monod} to ensure the existence of a non-zero fixed-point in $C$.

Let's prove the reverse implication. We show that $G$ has the translate property for $\RUC(G)$, which is equivalent to having the fixed-point property for cones for locally compact groups. Let $f\in \RUC(G)_+$ be a non-zero function and let $E = \Span_\R\left\{ gf: g\in G\right\}$ together with the supremum norm. Consider the cone $C = (E^*)_+$ in the locally convex vector space $\left( E^*, weak\text{-}* \right)$. Note that the adjoint action of $G$ on  $E^*$ is orbitally continuous for the weak-* topology, and $\left( E^*, weak\text{-}* \right)' = E$. Therefore, there is a non-zero positive $G$-invariant functional $\psi$ on $E$. By Proposition \ref{translate property implic che I(G) non vuoto}, this is sufficient to ensure that $f$ has the translate property.
\end{proof}

\subsection{The abstract translate property}

The subsection aims to distil the essence of the translate property for $\RUC(G)$. The strength of the translate property is that it doesn't depend on the ambient space. This fact will be helpful in situations where it is not possible to have control over the cobounded condition.

\begin{defn}
Let $G$ be a topological group. We say that $G$ has the \textbf{abstract translate property} if whenever $G$ has a continuous representation $\pi$ on a normed Riesz space $E$ by positive linear isometries, then for every non-zero $v\in E_+$
   \begin{align*}
        \sum_{j=1}^nt_jg_jv \geq 0 \quad \text{implies} \quad  \sum_{j=1}^n t_j \geq 0 \quad \text{ for every } t_1,...,t_n\in \R \text{ and } g_1,...,g_n \in G.
    \end{align*}
\end{defn}

Recall that the left-translation representation $\pi_L$ of $G$ on $\RUC(G)$ is the representation given by $\pi_L(g)f(x) = f(g^{-1} x)$, where $g,x\in G$ and $f\in \RUC(G)$, and the right-translation representation $\pi_R$ of $G$ on $\LUC(G)$ is the one given by $\pi_R(g)f(x) = f(xg)$, where $g,x\in G$ and $f\in \LUC(G)$.

\begin{lem}\label{translate property equal to abstract translate property}
Let $G$ be a topological group. Then $G$ has the abstract translate property if and only if $G$ has the translate property for $\mathcal{C}_{ru}^b(G)$.
\end{lem}
\begin{proof}
Suppose that $G$ has the abstract translate property. Then the left-translation representation of $G$ on the Banach space $\RUC(G)$ is continuous. Therefore, $G$ has the translate property for $\RUC(G)$.

Now, suppose that $G$ has the translate property for $\RUC(G)$ and let $\pi$ be a continuous representation of $G$ on a normed Riesz space $E$ by positive linear isometries. Let $v\in E$ be a non-zero positive vector and let $t_1,...,t_n\in \R$ and $g_1,...,g_n\in G$ such that $\sum_{j=1}^nt_j\pi(g_j)v \geq 0$. By Hahn-Banach Theorem, there is a positive linear functional $\lambda\in E'$ such that $\lambda(v) > 0$. Define the linear operator
    \begin{align*}
        T: E \longrightarrow \LUC(G), \quad w \longmapsto T(w),
    \end{align*}
where $T(w)(g) = \lambda(gw)$ for every $g\in G$. We should check that $T$ is well-defined. First, $T(w)$ is a bounded function for every $w\in E$, since $ \sup_{g\in G} |\lambda(gw)| \leq ||\lambda||_{op}||w||\text{ for every }g\in G$. Second, $T_v(w)$ is left-uniformly continuous. Indeed, let $\epsilon > 0$. As the representation of $G$ on $E$ is continuous, there is a neighborhood of the identity $A\subset G$ such that $ || aw - w || < \frac{\epsilon}{||\lambda||_{op}} $ for every $a\in A$.
Hence,
    \begin{align*}
       |\pi_R(a) T(w)(g)- T(w)(g)| 
        & = |\lambda(ga w)- \lambda(gw)|    
          = |\lambda(ga w-gw)|                 \\ 
        & \leq ||\lambda||_{op}||gaw-gw||     
        \leq ||\lambda||_{op}||aw-w|| < \epsilon
    \end{align*}
for every $g\in G$ and $a\in A$. This shows that $T(w)\in \LUC(G)$. Therefore, $T$ is well-defined. Notice that $T$ is equivariant with respect to the right-translation representation of $G$ on $\LUC(G)$. In fact,
    \begin{align*}
        T(aw)(g) = \lambda(gaw) = T(w)(ga) = \pi_R(a)T(w)(g)
    \end{align*}
for every $a,g\in G$ and $w\in E$.
Hence,
    \begin{align*}
        0 \leq T\left( \,\sum_{j=1}^nt_j\pi(g_j)v  \right) = \sum_{j=1}^n t_j\pi_R(g_j)T(v).
    \end{align*}
Nevertheless, saying that the representation $\pi_L$ of $G$ on $\RUC(G)$ has the translate property is equivalently saying that the representation $\pi_R$ of $G$ on $\LUC(G)$ has the translate property. Thus, $\sum_{j=1}^nt_j \geq 0$ as wished. 
\end{proof}

Note that the fixed-point property for cones implies the abstract translate property by Corollary \ref{corollary INI property implies TP}, and the two notions are equivalent for locally compact groups thanks to Theorem \ref{theorem TP imples measurably integral}.


\section{Hereditary Properties}\label{section hereditary properties}

We finally get to investigate the class of groups with the fixed-point property for cones. In the first part, we focus on the topological case, while in the second, we zoom on the locally compact one. Lastly, we give examples and non-examples of groups with the fixed-point property for cones. 

\subsection{The topological case}

The following proposition shows point b) of Theorem \ref{thm calss of top group with FPC}.

\begin{prop}\label{stable under continus homomorsphismus}
Let $G$ and $H$ be two topological groups and suppose that there is a surjective continuous group homomorphism $\phi: H \longrightarrow G$. Suppose that $H$ has the fixed-point property for cones, then $G$ has the fixed-point property for cones.
\end{prop}
\begin{proof}
We define the linear operator
    \begin{align*}
        T: \RUC(G) \longrightarrow \RUC(H), \quad h \longmapsto T (h) = h\circ \phi.
    \end{align*}
Then $T$ is a well-defined, non-zero and strictly positive. Moreover, it is equivariant. Indeed, let $\phi(g) \in G$ for some $g\in H$. Then
            \begin{align*}
                T\left(\phi(g)h\right)(x)   & = \big((\phi(g)h)\circ \phi\big)(x)
                                = (\phi(g)h)(\phi(x)) \\
                                & = h(\phi(g^{-1})\phi(x))
                                = h(\phi(g^{-1}x)) = gT (h)(x)
            \end{align*}
        for every $h\in \mathcal{C}_{ru}^b(G)$ and every $x\in H$.
        
    Now the function $T(f)\in \RUC(H)$ is non-zero and positive for every non-zero $f\in \RUC(G)_+$, therefore the application $T$ maps $\RUC(G,f)$ to a subspace of $\RUC(H,T(f))$. On this last space, there is an invariant normalized integral $\I$ by hypothesis. The composition of $\I$ with $T$ gives an invariant normalized integral for $\RUC(G,f)$. Consequently, $G$ has the invariant normalized integral property for $\RUC(G)$ and hence the fixed-point property for cones by Theorem \ref{FPC equivalente integral}.
\end{proof}

Direct consequences are:

\begin{cor}
Suppose that $G$ admits a topology $\tau_1$ for which $(G,\tau_1)$ has the fixed-point property for cones. Then $(G,\tau_2)$ has the fixed-point point property for cones for every weaker topology $\tau_2$.
\end{cor}
\begin{proof}
This is only because that the identity map $\Id:(G,\tau_1) \longrightarrow (G,\tau_2)$ is a continuous isomorphism of groups, since the topology $\tau_2$ is weaker than the topology $\tau_1$.
\end{proof}

\begin{cor}
If $G$ has the fixed-point property for cones as an abstract group, i.e., when equipped with the discrete topology, then $G$ has the fixed-point property for cones for every admissible group topology, i.e., every topology which makes $G$ a topological group.
\end{cor}

Moreover, the fixed-point property for cones is preserved under taking quotients, since the quotient map $p: G \longrightarrow \faktor{G}{N}$, where $G$ is a topological group and $N$ a normal subgroup, is a continuous and surjective homomorphism of groups (\cite[Theorem $(5.16)$]{abstractharmonic}). 

\medskip
The proof of point a) of Theorem \ref{thm calss of top group with FPC} is provided in the following proposition.

\begin{prop}
Let $G$ be a topological group which has the fixed-point property for cones and let $H\leq G$ be an open subgroup. Then $H$ has the fixed-point property for cones. 
\end{prop}

\begin{proof}
We show that $H$ has the invariant normalized integral property for $\RUC(H)$. To this aim, let $K$ be a set of representatives for the right $H$-cosets and define the linear operator
    \begin{align*}
        T: \RUC(H) \longrightarrow \RUC(G), \quad h \longmapsto T (h)(g) = h(x),
    \end{align*}
where $g = xk$ for a $k\in K$. Then $T$ is well-defined, strictly positive and equivariant. As done above, this implies that $H$ has the invariant normalized integral property for $\RUC(H)$ as $G$ has the invariant normalized integral property for $\RUC(G)$. Therefore, $H$ has fixed-point property for cones by Theorem \ref{FPC equivalente integral}. 
\end{proof}

The following result completes the proof of Theorem \ref{thm calss of top group with FPC}.

\begin{prop}
Let $G$ be a topological group and let $F \unlhd G$ be a finite normal subgroup. Suppose that the group $\faktor{G}{F}$ has the fixed-point property for cones. Then $G$ has the fixed-point property for cones.
\end{prop}
\begin{proof}
Suppose that $G$ has representation on a non-empty proper weakly complete cone $C$ in a locally convex vector space $E$ which is of cobounded type and locally bounded right-uniformly continuous. We have to show that $G$ fixes a non-zero vector in $C$. Consider the closed subspace $E^F$ of all vectors of $E$ which are $F$-invariant and the natural representation of $\faktor{G}{F}$ on it.
Then:
    \begin{itemize}
        \item[-] the proper cone $C^F = C \cap E^F$ is weakly complete as the intersection of complete sets is complete and non-empty because $\sum_{g\in F} gc \in C^F$ for every $c\in C$;
        
        \item[-] the representation of $\faktor{G}{F}$ on $E^F$ is locally bounded right-uniformly continuous as witness by the vector $x_F = \sum_{g\in F} gx_0$, where $x_0$ is the vector which gives the locally bounded right-uniformly condition for the representation of $G$ on $E$. Note that $x_F \neq 0$ as the cone $C$ is proper;
        
        \item[-] the representation of $\faktor{G}{F}$ on $E^F$ is also of cobounded type because the map $E' \longrightarrow (E^F)'$ is equivariant, positive and onto by Hahn-Banach Theorem.
    \end{itemize}
Therefore, the fixed-point property for cones of $\faktor{G}{F}$ gives a non-zero $\faktor{G}{F}\,$-fixed-point $c_0\in C^F$. The vector $c_0$ is also a non-zero $G$-fixed-point. 
\end{proof}

To conclude the discussion about topological groups, we point out that the fixed-point property for cones is invariant under passing to finite-index subgroups and passes to dense subgroups. 

\begin{prop}\label{finite-index subgroups to ambient groups FPC}
Let $G$ be a topological group and let $H < G$ be a topological subgroup of finite-index. If $H$ has the fixed-point property for cones, then $G$ has the fixed-point property for cones.
\end{prop}
\begin{proof}
Consider a locally bounded right-uniformly continuous representation of $G$ on a non-empty weakly complete proper cone $C$ which is of cobounded type. To apply the fixed-point property for cones of $H$ on $C$, we only have to justify that the representation of $H$ on $C$ is of cobounded type.  Let $\lambda\in E'$ be the $G$-dominating element of $E'$ and fix $R\subset G$ a set of representatives for the right $H$-cosets in $G$. We claim that $\lambda_H = \sum_{r\in R}r\lambda$ is a $H$-dominating element for $E'$. Indeed, any finite sum $\sum_{j=1}^n g_j\lambda$, where $g_1,...,g_n\in G$, can be written as $\sum_{j=1}^n h_jr_j\lambda$ for $h_1,...,h_n\in H$ and $r_1,...,r_n\in R$. This last sum is dominated by $\sum_{j=1}^n h_j \lambda_H$ because $r_j\lambda \leq \lambda_H$ for every $r_j\in R$. Therefore, there is a non-zero $H$-fixed-point $x_H$ in $C$. The vector $x = \sum_{r\in R} r^{-1}x_H$ is a $G$-fixed-point in $C$ and $x \neq 0$, since the cone $C$ is proper. 
\end{proof}

\begin{prop}\label{dense subgroup of topological with fixed-point property}
Let $G$ be a topological group and let $D < G$ be a dense topological subgroup. If $G$ has the fixed-point property for cones, then $D$ has the fixed-point property for cones. 
\end{prop}
\begin{proof}
Consider the extension operator 
    \begin{align*}
        \texttt{ext}: \RUC(D) \longrightarrow \RUC(G), \quad f \longmapsto \texttt{ext}(f),
    \end{align*}
which gives to every bounded right-uniformly continuous map $f$ on $D$ its unique extension to $G$, see \cite[Chap.II §3 No.6 Théorème 2]{bourbakitopgen}. Then $\texttt{ext}$ is a positive $D$-equivariant  linear operator. Take a non-zero $f\in \RUC(D)_+$. We want to show that there is an invariant normalized integral on $\RUC(D,f)$. As $G$ has the fixed-point property for cones, there is a $G$-invariant normalized integral $\I$ on $\RUC(G, \texttt{ext}(f))$. Thus, $\overline{\I}:= \I \circ \texttt{ext}$ is a $D$-invariant normalized integral on $\RUC(D,f)$.
\end{proof}

We could not prove or disprove that the fixed-point property for cones is invariant under passing to dense subgroups. The central problem in finding a proof is the absence of continuity of the group representations when working with $p_d$-norms. 

\medskip
In the setting of topological groups, the following result is a priori weaker than saying that the fixed-point property for cones is preserved by directed union. 
\medskip
Recall that having the translate property is equivalent to having the abstract translate property by Lemma \ref{translate property equal to abstract translate property}.

\begin{prop}\label{directed unin of topological subgroup}
Let $G$ be a topological group. Suppose that $G$ is the directed union of a family $(H_\alpha)_\alpha$ of topological subgroups. If every $H_\alpha$ has the abstract translate property, then $G$ has the translate property for $\RUC(G)$.
\end{prop}
\begin{proof}
Let $f\in \RUC(G)_+$ be a non-zero function and let $t_1,...,t_n\in \R$ and $g_1,...,g_n\in G$ such that $\sum_{j=1}^nt_jg_jf \geq 0$. We have to show that $\sum_{j=1}^nt_j \geq  0$. Take $\alpha$ such that $g_1,...,g_n\in H_\alpha$. Then the abstract translate property of $H_\alpha$ ensures that  $\sum_{j=1}^nt_j \geq  0$ as wished.
\end{proof}

\subsection{The locally compact case}
Let's now look at the case of locally compact groups. It is possible to study deeper the class of locally compact groups with the fixed-point property for cones thanks to the results of the Sections \ref{Chapter the locally compact case} and \ref{section equivalent fixed-point properties}.

\medskip
We start by showing Theorem \ref{Closed group have the FPC}. To this end, recall that if $G$ is a locally compact group and $H<G$ a closed subgroup, then a Bruhat function $\beta$ for $H$ is a positive continuous function on $G$ such that $\supp(\beta|_{KH})$ is compact for every $K$ compact subset of $G$ and such that $\int_{H}\beta (gh)dm_H(h)=1$ for every $g\in G$. Note that $\beta$ depends of the choice of the Haar measure $m_H$ of $H$. A proof that every closed subgroup of a locally compact group admits a Bruhat function $\beta$ can be found in \cite[Proposition $1.2.6$]{runde}.

\begin{proof}[Proof of Theorem \ref{Closed group have the FPC}]
We want to show that $H$ has the translate property for $L^\infty(H)$. Therefore, let $\beta$ be a Bruhat function of $H$ and define the operator
    \begin{align*}
        T: L^\infty (H) \longrightarrow L^\infty(G), \quad h \longmapsto T(h)(g) = \int_H h(x)\beta(g^{-1}x)dm_H(x).
    \end{align*}
Note that $T$ is linear and well-defined as $T(h)\in \CB(G)$ for every $h\in L^\infty(H)$, see \cite[Proposition $(1.12)$]{paterson}. Moreover, $T$ is strictly positive as $\beta$ is positive and its support intersects that of $f$. Finally, $T$ is equivariant. Indeed,
    \begin{align*}
        T(ah)(g) & = \int_H (ah)(x)\beta(g^{-1}x)dm_H(x)  
                   = \int_H h(a^{-1}x)\beta(g^{-1}x)dm_H(x)  \\
                 & = \int_H h(y)\beta(g^{-1}ay)dm_H(y)  
                   = \int_H h(y)\beta((a^{-1}g)^{-1}y)dm_H(y)  \\
                 & = T(h)(a^{-1}g) = aT(h)(g)
    \end{align*}
for every $a\in G$ and $h\in L^{\infty}(H)$.
Now let $f\in L^\infty(H)$ be a non-zero positive function and let $h_1,...,h_n \in H$ and $t_1,...,t_n \in \R$ be such that $\sum_{j=1}^n t_jh_jf \geq 0$. Therefore,
    \begin{align*}
        0 \leq T\Big( \sum_{j=1}^n t_jh_jf \Big) = \sum_{j=1}^n t_jh_jT(f).
    \end{align*}
As $G$ has the fixed-point property for cones, then $G$ also has the translate property for $\CB(G)$ by Theorem \ref{Equivalence of TP on different spaces}. Applying it to the non-zero positive function $T(f)$, we have that $\sum_{j=1}^n t_j \geq 0$. This proves that $H$ has the translate property for  $L^\infty(H)$ and hence the fixed-point property for cones by Theorems \ref{TP on E equivalent to NII on E} and \ref{Equivalence of integral on different spaces}.
\end{proof}

We continue illustrating the proof of Theorem \ref{ thm in loc cpt for dense/union/lattices}. The proof of point a) is a consequence of Proposition \ref{directed unin of topological subgroup} because the translate property for $\RUC(G)$ is equivalent to the fixed-point property for cones for locally compact groups (Theorems \ref{FPC equivalente integral} and \ref{theorem TP imples measurably integral}). In particular, the fixed-point property for cones for locally compact groups is a local property, i.e., a locally compact group has the fixed-point property for cones if and only if every compactly generated subgroups has it. The proof of point b) is given in the following proposition. 

\begin{prop}
Let $G$ be a locally compact group and let $D < G$ be a dense subgroup. Then $D$ has the fixed-point property for cones if and only if $G$ has the fixed-point property for cones.
\end{prop}
\begin{proof}
The $\textit{if}$ direction is given by Proposition \ref{dense subgroup of topological with fixed-point property}. 

So let's prove the \textit{only if} direction. Suppose that $G$ has a continuous representation on a normed Riesz space $E$ by positive linear isometries. Let $v\in E_+$ be a non-zero vector and take $g_1,...,g_n\in G$ and $t_1,...,t_n\in \R$ such that $\sum_{j=1}^n t_jg_jv \geq 0$. Since $D$ is dense in $G$, for every $j$ there is a net $(d_\alpha^{(j)})_\alpha \subset D$ such that $\lim_\alpha d_\alpha^{(j)} = g_j$. Moreover, $\lim_\alpha d_\alpha^{(j)}v = g_jv $ for every $j$ as the representation of $G$ on $E$ is continuous. Therefore, 
    \begin{align*}
       0 \leq  \sum_{j=1}^n t_jg_jv  = \sum_{j=1}^n t_j\lim_\alpha d_\alpha^{(j)}v = \lim_\alpha \sum_{j=1}^n t_jd_\alpha^{(j)}v.
    \end{align*}
The abstract translate property of $H$ implies that $\sum_{j=1}^n t_j \geq 0$. Hence, $G$ also has the abstract translate property. However, for locally compact groups the abstract translate property is equivalent to the fixed-point property for cones by Lemma \ref{translate property equal to abstract translate property} and Theorem \ref{theorem TP imples measurably integral}. Thus, $G$ has the fixed-point property for cones. 
\end{proof}

\subsection{Well-behaved group extensions} As previously said, the fixed-point property for cones is not preserved by group extensions in general. However, some nice extensions preserve it as theorem \ref{cartesian product with subexponential growth} witness.

\medskip
We recall that a locally compact group $G$ is of \textbf{subexponential growth} if $\lim_n m_G(C^n)^{\frac{1}{n}}=1$ for every compact neighborhood of the identity $C\subset G$.

\medskip
The proof of point a) of Theorem \ref{cartesian product with subexponential growth} is actually the same as the one given in \cite[Theorem $8-(3)$]{monod}. Note that instead of $\ell^\infty(G)$, one should use $\RUC(G)$ and the lemma employed in Monod's proof (\cite[Lemma 30]{monod}) is actually true for locally compact groups. In fact, it was proved by Jenkins in \cite[Lemma 1]{jenkinsexponential}.

\medskip
In particular, point a) of Theorem \ref{cartesian product with subexponential growth} implies that every locally compact group with subexponential growth has the fixed-point property for cones. This leads us to the following two corollaries. 

\begin{cor}
Virtually nilpotent locally compact groups have the fixed-point property for cones.
\end{cor}
\begin{proof}
Nilpotent locally compact groups have subexponential growth by \cite[(6.18)]{paterson}. Moreover, the fixed-point property passes through finite-index subgroups by Proposition \ref{finite-index subgroups to ambient groups FPC}.
\end{proof}

We say that a locally compact group $G$ is a \textbf{topologically finite conjugancy classes group}, or $G$ is a \textbf{$\overline{FC}$-topologically group}, if the closure of each of its conjugancy classes is compact. These types of groups have subsexponential growth, see \cite{palmer}. Hence,

\begin{cor}
Let $G$ be a locally compact $\overline{FC}$-topologically group. Then $G$ has the fixed-point property for cones.
\end{cor}

\medskip
We continue by demonstrating point b) of Theorem \ref{cartesian product with subexponential growth}.

\begin{proof}[Proof of Theorem \ref{cartesian product with subexponential growth} - b)]
Let $G$ be a locally compact group which is an extension of a compact group $C$ by a group with the fixed-point property for cones $Q$. We want to show that $G$ has the translate property for $\RUC(G)$.

Set $E = \RUC(G)$ and define the strictly positive linear operator given by the Bochner integral
    \begin{align*}
        T: E \longrightarrow E^C, \quad f \longmapsto T(f) = \int_C cf \; dm_C(c),
    \end{align*}
where $m_C$ is the normalized Haar measure of $C$. Then $T$ is well-defined as
    \begin{align*}
        xT(f) = \int_C xcf \; dm_C(c) = \int_C yf \; dm_C(x^{-1}y) = \int_C yf \; dm_C(y) =  T(f)
    \end{align*}
for every $f\in E$ and $x\in C$.
Moreover, $T$ is equivariant with respect to the natural representation of the quotient group $Q$ on $E^C$. Indeed, for every $g\in G$ and every $c\in C$, the element $g^{-1}cg$ is in $C$ by normality. Therefore, for every $c\in C $ there is $c'\in C$ such that $cg = gc'$. Thus,
    \begin{align*}
        T(gv) = \int_C cgv \; dm_C(c) = \int_C gc'v \; dm_C(gcg^{-1}) = g \int_C c'v\; dm_C(c) = gT(v)
    \end{align*}
for every $v\in E$. Now let $f\in E_+$ be a non-zero function and let $t_1,...,t_n\in \R$ and $g_1,...,g_n\in G$ such that $\sum_{j=1}^n t_jg_jf \geq 0$. Then
    \begin{align*}
        0 \leq T \left( \sum_{j=1}^n t_jg_jf  \right) = \sum_{j=1}^n t_jq_jT(f)
    \end{align*}
for some $q_1,...,q_n\in Q$. Applying the abstract translate property of $Q$ to the function $T(f)$, it is possible to conclude that $\sum_{j=1}^n t_j \geq 0$. This shows that $G$ has the translate property for $\RUC(G)$ and hence the fixed-point property for cones.  
\end{proof}

\begin{scholium}
There is a direct way to show that a locally compact group with subexponential growth has the fixed-point property for cones. The strategy was pointed out by Paterson in \cite[(6.42)(i)]{paterson} which claimed that a locally compact group with subexponential growth has the \textit{generalized translate property} and so, consequently, the fixed-point property for cones (see Scholium \ref{schoulium 1}).
\end{scholium}

\subsection{Obstruction to the fixed-point property for cones}\label{Obstruction to the fixed point property on convex cones} We discussed groups with the fixed-point property for cones, but let's now discuss the groups that don't have it. To begin, look at the following obstruction. 

\medskip
We recall that a \textbf{uniformly discrete free subsemigroup in two generators} of a group $G$ is a subsemigroup $\text{T}_2$ generated by two elements $a,b\in G$ such that there is a neighborhood $W$ of the identity with the property that $sW\cap tW = \emptyset$ for every $s,t\in \text{T}_2$, $s\neq t$.

\begin{prop}\label{T_2 obstruction to be (FPC)}
Let $G$ be a locally compact group which contains a uniformly discrete free subsemigroup in two generators $\text{T}_2$. Then $G$ doesn't have the fixed-point property for cones.
\end{prop}
\begin{proof}
Suppose it is not the case. Therefore, there is a neighborhood $W$ of the identity such that $sW\cap tW = \emptyset$ for every $s,t\in \text{T}_2$. Define the open subset $U:= \text{T}_2\cdot W$ and note that $U$ has Haar measure bigger than zero. Then $\1_U$ is different to the zero function. By Theorem \ref{Equivalence of integral on different spaces}, there is an invariant normalized integral $\I$ on $L^\infty(G,\1_U)$. But now $aU\cap bU =\emptyset$. Therefore, the function $\phi := \1_U -\1_{aU}-\1_{bU}$ is non-zero and positive.

Note that $\phi\in L^{\infty}(G,\1_U)$. Therefore,
    \begin{align*}
        0 \geq \I(\phi) = \I(\1_U) - \I(\1_{aU})- \I(\1_{bU}) = -1,
    \end{align*}
which is a contradiction.
\end{proof}

In particular, this last definition shows that not every extension of groups with the fixed-point property for cones has the fixed-point property for cones. The easiest example is given by the discrete group of affine transformations of the line $\R \rtimes \R^*$. In fact, it contains the free subsemigroup $\text{T}_2$ generated by the elements $(0,2)$ and $(1,1)$.

\medskip
Thanks to this obstruction, it is possible to show that there are exceptional cases where the fixed-point property for cones is equivalent to being of subexponential growth. 

\medskip
Let $G$ be a locally compact group and write $G_e$ for the connected component of the neutral element. Then $G$ is said \textbf{almost connected} if the topological quotient group $\faktor{G}{G_e}$ is compact. 

\begin{prop}\label{connectec FPC equivalent to subexp}
Let $G$ be a locally compact group. 
    \begin{itemize}
        \item[a)] Suppose that $G$ is connected. Then $G$ has the fixed-point property for cones if and only if $G$ has subexponential growth.
        \item[b)] Suppose that $G$ is compactly generated and almost connected. Then $G$ has the fixed-point property for cones if and only if $G$ has subexponential growth.
    \end{itemize}
\end{prop}
\begin{proof}
The \textit{if} direction is true in both cases by Proposition \ref{cartesian product with subexponential growth}.

For the \textit{only if} direction of point a), suppose that $G$ is not of subexponential growth. Then there is a uniformly discrete subsemigroup $\text{T}_2\subset G$ in two generators by \cite[Proposition 6.39]{paterson}. By Proposition \ref{T_2 obstruction to be (FPC)}, $G$ can not have the fixed-point property for cones.

For the \textit{only if} direction of  point b), suppose that $G$ has the fixed-point property for cones and consider the exact sequence given by
    \begin{align*}
        \{e\} \longrightarrow G_e \longrightarrow G \longrightarrow \faktor{G}{G_e} \longrightarrow \{ e\},
    \end{align*}
where $G_e$ is the connected component of the identity. As $G_e$ is closed and $\faktor{G}{G_e}$ is compact, the groups $G$ and $G_e$ have the same growth by \cite[Theorem 4.1]{guivarch}. But now $G_e$ has the fixed-point property for cones by Proposition \ref{Closed group have the FPC}, and therefore it also has subexponential growth by point a). Thus, we can conclude that $G$ has subexponential growth.
\end{proof}


\section{Fixing Radon Measures}\label{section fixing radon measures}

The section aims to give a couple of applications of the fixed-point property for cones. In particular, we will see that the fixed-point property for cones is the correct property to work with problems where a non-zero invariant Radon measure is required.

\subsection{Cocompact actions and fixed Radon measures}

Let $G$ be a locally compact group that acts on a locally compact space $X$. 

\medskip
Recall that the action of $G$ on $X$ is called \textbf{cocompact} if there is a compact subset $K$ of $X$ such that the image of $K$ under the action of $G$ covers $X$.

\begin{prop}\label{prop cocompact action and cobounded}
The group $G$ acts cocompactly on $X$ if and only if the space $\mathcal{C}_{00}(X)$ admits a $G$-dominating element.
\end{prop}
\begin{proof}
Suppose that the action of $G$ on $X$ is cocompact. Then there is a compact subset $K\subset X$ such that $X = \bigcup_{g \in G}gK$. Take a relatively compact neighborhood $U_k$ of $k$ for every $k\in K$. Then $K \subset \bigcup_{k\in K}U_k$. As $K$ is compact, there are $k_1,...,k_n \in K$ such that $K\subset \bigcup_{j=1}^nU_{k_j}$. Define $U = \bigcup_{j=1}^nU_{k_j} $. Because every $U_{k_j}$ is relatively compact then so is $U$. Thus, there is an open set $V\subset X$ such that $\overline{U}\subset V$. By Uryshon Lemma (\cite[Lemma $2.12$]{rudinanal}), there is  $\psi\in \mathcal{C}_{00}(X)_+$ such that $\psi = 1$ on $\overline{U}$ and $\psi = 0 $ on $X \setminus V$. We claim that $\psi$ is a $G$-dominating element in $\mathcal{C}_{00}(X)$. Indeed, let $\phi\in \mathcal{C}_{00}(X)$ and let $K'= \supp(\phi)$. As the action is cocompact, there are $g_1,...,g_n\in G$ such that 
    \begin{align*}
        K'\subset \bigcup_{j=1}^ng_jK \subset \bigcup_{j=1}^ng_j\overline{U}.
    \end{align*}
This implies that
    \begin{align*}
        |\phi| \leq ||\phi||_{\infty} \1_{supp(\phi)}
               \leq \sum_{j=1}^n ||\phi||_{\infty} \1_{g_j\overline{U}}
               = \sum_{j=1}^n ||\phi||_{\infty} g_j\1_{\overline{U}}
               \leq \sum_{j=1}^n ||\phi||_{\infty} g_j \psi
    \end{align*}
as wanted.

Let now suppose that $\mathcal{C}_{00}(X)$ admits a $G$-dominating element $\psi$, and we want to show that the action of $G$ on $X$ is cocompact. To this aim, it suffices to show that there is a compact set $K\subset G$ such that for every $x\in X$ there is $g\in G$ with $x\in gK$. Define $K = \supp(\psi)$. Take $x\in X$ and $K'\subset X$ a compact neighborhood of $x$. By Urysohn Lemma, there exists $\phi \in \mathcal{C}_{00}(X)$ with $K' \subset \supp(\phi)$. Thus, there are $g_1,...,g_n\in G$ such that $|\phi|\leq \sum_{j=1}^ng_j\psi$. But this implies that $K' \subset \bigcup_{j=1}^ng_jK$, and so we are done.
\end{proof}

Moreover, it is not difficult to see that it is possible to chose the dominating element of $\mathcal{C}_{00}(X)$ to be support-dominating.

\begin{thm}\label{invariant Radon measure for cocompact -FPC}
Suppose that $G$ is a locally compact group with the fixed-point property for cones. Then each orbitally continuous cocompact action of $G$ on any locally compact topological space $X$ has a non-zero invariant Radon measure.
\end{thm}

\begin{proof}
Recall that each Radon measure on a locally compact space $X$ can be seen as a positive functional defined on the Riesz space $\mathcal{C}_{00}(X)$ (\cite[Chap. III §1 No.5 Théorème 1]{bourbakihaar}). Therefore, consider the cone $C = \mathcal{C}_{00}(X)_+^*$ in the locally convex vector space $E = (\mathcal{C}_{00}(X)^*, $ weak-*$)$. This cone is weakly complete as it is weak-* closed in the complete space $E$, and the adjoint representation of $G$ on $C$ is of cobounded type because $E' = \mathcal{C}_{00}(X)$, and this last space admits $G$-dominating element by Proposition \ref{prop cocompact action and cobounded}. Moreover, the representation of $G$ on $C$ is continuous. Indeed, the action of $G$ on $X$ is also jointly continuous by Ellis Theorem (\cite[Lemma C.1.9]{max}). This implies that the induced representation of $G$ on $\mathcal{C}_{00}(X)$ is $||\cdot||_\infty$-continuous, and so the adjoint representation of $G$ on $E$ is weak-* continuous. Therefore, point b) of Lemma \ref{lemma continuous fixed-point property for cones} provides a non-zero fixed-point in $C$ which is nothing but a non-zero invariant Radon measure. 
\end{proof}

\subsection{Unimodularity}

Let $G$ be a locally compact group and write $\Delta_G$ for its unimodular function.

\begin{cor}\label{unimodularity of FPC-SUPER groups}
Let $G$ be a locally compact group with the fixed-point property for cones, then $G$ is unimodular.
\end{cor}

\begin{proof}
Suppose it is not the case. This means that there is $g\in G$ such that $\Delta_G (g) = c \neq 1$. Let $H:= \langle g \rangle$ be the group generated by $g$. We claim that $H\cong \Z$, and that $H$ is closed and discrete as a subgroup of $G$. This is because of the following facts.
\begin{itemize}
    \item[-] \textit{The element $g$ has infinite order.} Suppose it is not the case. Then there is $n\in \N$ such that $g^n=e$. This implies that 
    \begin{align*}
        \Delta_G(g)^n = \Delta_G(g^n) = \Delta_G (e) = 1 \quad \iff  \quad \Delta_G (g) = 1.
    \end{align*}
    But this is a contradiction.
    
    \item[-] \textit{$H$ is closed.} Suppose it is not the case. Then there is a net $(g_\alpha)_\alpha \subset H$ such that $\lim_\alpha g_\alpha = x$ and $x\notin H$. Note that $(g_\alpha)_\alpha = (g^{n_\alpha})_\alpha$, where $(n_\alpha)_\alpha$ is a net in $\N$ converging to infinity. Then 
        \begin{align*}
            \Delta_G(x) =\Delta_G( \lim_{\alpha} g^{n_\alpha})
                        = \lim_{\alpha} \Delta_G(g)^{n_\alpha}
                        = \lim_{\alpha} c^{n_\alpha} = 
                        \begin{cases}
                            0       & \quad \text{if } c<1\\
                            \infty  & \quad \text{if } c>1.
                        \end{cases}
        \end{align*}
    But this is a contradiction with the fact that $im(\Delta_G)$ is a subgroup of $\R^*$.
    
    \item[-] \textit{$H$ is discrete.} Suppose it is not the case. Then there is a net $(g_\alpha)_\alpha \subset H$ such that $\lim_\alpha g_\alpha = e$. Employing the same strategy as before, we can find a contradiction.
\end{itemize}
Therefore, we can conclude that $H$ is isomorphic to $\Z$. Now let $G$ act on the locally compact space $\faktor{G}{H}$, and note that this action is continuous and cocompact. By Theorem \ref{invariant Radon measure for cocompact -FPC}, there is a non-zero invariant Radon measure on  $\faktor{G}{H}$. But this is possible if and only if the restriction of $\Delta_G$ on the subgroup $H$ is equal to $\Delta_H$ as explained in \cite[Chap. II §2 No.6 Corollaire 2]{bourbakiintegral}. As $H$ is discrete, $\Delta_H = 1$. Thus,
    \begin{align*}
        1\neq \Delta_G(g) = \Delta_H(g) = 1,
    \end{align*}
which is a contradiction. We can conclude that $G$ is unimodular.
\end{proof}

\begin{cor}
Closed subgroups of groups with the fixed-point property for cones are unimodular.
\end{cor}

\begin{proof}
Closed subgroups of a group with the fixed-point property for cones have the fixed-point property for cones by Theorem \ref{Closed group have the FPC}. Thus, they are unimodular.
\end{proof}

\begin{rem}
In general, unimodularity doesn't pass to closed subgroups. An easy example is given by the locally compact Lie group $\GL_2(\R)$ which is unimodular, but it contains the non-unimodular $(ax+b)$-group.
\end{rem}

\subsection{Fixing Radon measures on the line}

We want to apply the close relationship between the fixed-point property for cones and non-zero invariant Radon measures to a dynamical problem. Namely, when the natural action of a subgroup of $\Homeo(\R)$, the group of order-preserving homeomorphisms of the line, on $\R$ fixes a non-zero Radon measure. This problem was primarily studied by Plante (\cite{plante} and \cite{plantefoliations}), who discovered that finitely generated virtually nilpotent subgroups of $\Homeo(\R)$ always fix a non-zero Radon measure on the line (an exposition of Plante's work can be found in \cite[Subsection 2.2.5]{navas}). This result is due to the fact that finitely generated virtually nilpotent groups are of subexponential growth, so it is not surprising that we can generalise this theorem to the class of groups with the fixed-point property for cones. Interestingly, the proof for groups with the fixed-point property for cones is much more natural and less technical than the one for virtually nilpotent groups.

\medskip
Recall that a subgroup $G$ of $\Homeo(\R)$ is said to be \textbf{boundedly generated} if there is a symmetric set of generators $S$ of $G$ and a point $x_0\in \R$ such that the set $\{ sx_0 : s\in S \}$ is a bounded subset of the line.

\begin{thm}\label{boundely generated and fixed-radon measure}
Let $G$ be a boundedly generated subgroup of $\Homeo(\R)$ with the fixed-point property for cones. Then there is a non-zero $G$-invariant Radon measure on $\R$.
\end{thm}

\begin{proof}
Let $S\subset G$ be a symmetric set that generates boundedly $G$. 

Clearly, if the action of $G$ on $\R$ has a global fixed-point, then a Dirac mass on this point is a non-zero invariant Radon measure. 

Therefore, suppose that the action of $G$ on $\R$ has no global fixed-points. We claim that the action is of cocompact-type. Indeed, let $x_0\in \R$ be the point which witnesses the fact that $S$ generates $G$ boundedly and let $x_1 := \sup_{s\in S} s_jx_0$. Then the set $ I = [x_0,x_1]$ is compact because $S$ generates $G$ boundedly by hypothesis. If we can show that every orbit of $G$ intersects the interval $I$, then we have cocompactness. Let $x\in \R$ and note that the orbit $Gx$ is unbounded in the two directions, otherwise its supremum and its infimum would be global fixed-points. Thus, we can chose $x_0',x_1'$ in $Gx$ such that $x_0' < x_0 < x_1 <x_1'$. Let $g = s_{j_n}\cdots s_{j_1} \in G$ be such that $gx_0' = x_1'$ and let $m \in \{ 1,...,n-1\}$ be the largest index for which the inequality $s_{j_{m}} \cdots s_{j_1}x_0' < x_0$ holds. Then $s_{j_{m+1}}s_{j_m}\cdots s_{j_1}x_0'$ is in the orbit of $x$, it is greater than or equal to $x_0$ and it is smaller or equal to $x_1$ by definition. Therefore, the orbit of $x$ intersects $I$. Now Theorem \ref{invariant Radon measure for cocompact -FPC} provides a non-zero invariant Radon measure as wanted. 
\end{proof}

In particular, Theorem \ref{boundely generated and fixed-radon measure} generalizes \cite[Theorem $(5.4)$]{plantefoliations}.

\medskip

\begin{rem}
The assumption that $G$ is boundedly generated cannot be dropped. In fact, there are even examples of abelian subgroups of $\Homeo(\R)$ which have no non-zero invariant Radon measure on the line (see \cite[Section $5$]{plante}).
\end{rem}

\begin{scholium}
The theory developed in this last subsection could also be developed using the slightly more general concept of supramenability instead of the fixed-point property for cones. This is possible thanks to \cite[Proposition 2.7]{superamenablegroups}.
\end{scholium}



\begin{thebibliography}{9}

\bibitem[AB99]{infinitedim} C.D. Aliprantis \& K.C. Border, \emph{Infinite dimensional analysis : a hitchhiker's guide}, Springer, Berlin-Heidelberg, 1999. 

\bibitem[AT07]{conesandduality} C.D. Aliprantis \& R. Tourky, \emph{Cones and Duality}, Graduate Studies in Mathematics, Vol. 84, 2007.

\bibitem[Bou59]{bourbakiintegral} Nicolas Bourbaki, \emph{Éléments de mathématique. Livre VI: Intégration}, Hermann, 1959.

\bibitem[Bou63]{bourbakihaar} Nicolas Bourbaki, \emph{Éléments de mathématique. Livre VII: Intégration}, Hermann, 1963.

\bibitem[Bou71]{bourbakitopgen} Nicolas Bourbaki, \emph{Éléments de mathématique. Livre III: Topologie générale,}, Chapitres 1 à 4, Hermann, 1971.

\bibitem[Bou81]{bourbakiespace} Nicolas Bourbaki, \emph{Éléments de mathématique. Livre VI: Espaces vectoriels topologiques}, Chapitres 1 à 5, Springer, 1981.

\bibitem[C13]{measurecohn} D. L. Cohn, \emph{Measure Theory}, Springer Science \& Business Media, 2013.

\bibitem[DE09]{principle} A. Dietmar \& S. Echterhoff, \emph{Principle of harmonic analysis}, Springer, New York, 2009.

\bibitem[D61]{day} Mahlon M. Day, \emph{Fixed-point theorems for compact convex sets}, Illinois Journal of Mathematics, Vol. 5, Issue 4, pp. 585-590, 1961.

\bibitem[G17]{max} Maxime Gheysens, \emph{Representing groups again all odds}, PhD's thesis, EPFL, 2017.

\bibitem[G69]{greenleaf} Frederick P. Greenleaf, \emph{Invariant Means on Topological Groups and their Applications}, Van Nostrand-Reinhold Mathematical Studies, 1969. 

\bibitem[G73]{guivarch} Yves Guivarc'h, \emph{Croissance polynomiale et périodes des fonctions harmoniques}, Bulletin de la S.M.F., Vol. 101, pp. 333-379, 1973.

\bibitem[HR63]{abstractharmonic} E. Hewitt, K. A. Ross, \emph{Abstract Harmonic Analysis I}, 2 Edition, Grundlehren der mathematisch Wissenschaften, Vol. 115, Springer, New York, 1963.

\bibitem[HR70]{abstractharmonic2} E. Hewitt, K. A. Ross, \emph{Abstract Harmonic Analysis II}, 2 Edition, Grundlehren der mathematisch Wissenschaften, Vol. 152, Springer, New York, 1970.

\bibitem[J74]{jenkisfolner} J.W. Jenkins, \emph{F\/olner's condition for exponentially bounded groups}, Mathematica Scandinavica, Vol. 35, No. 2, pp. 165-174, 23 April 1975

\bibitem[J76]{jenkinsexponential} J.W. Jenkins, \emph{A Fixed Point Theorem for Exponentially Bounded Groups}, Journal of Functional Anylsis, Vol. 22, No. 4, pp. 346-353, August 1976.

\bibitem[J80]{jenkinsfixedpoint} J.W. Jenkins, \emph{On Groups Action With Nonzero Fixed Point}, Pacific Journal of Mathematics, Vol. 91, No. 2, 1980.

\bibitem[K14]{kellerhals} Julian Kellerhals, \emph{Supramenable groups}, PhD thesis, EPFL, Lausanne, 2014.

\bibitem[KMR13]{superamenablegroups} J. Kellerhals, N. Monod \& M. R{\o}rdam, \emph{Non-supramenable groups acting on locally compact spaces}, Doc. Math., Vol. 18, pp. 1597–1626, 2013.

\bibitem[L90]{lau} A. to-Ming Lau, \emph{Amenability for Semigroups}, The Analytical and Topological Theory of Semigroups: Trends and Developments, de Gruyter Expositions in Mathematics, Vol. 1, Berlin, 1990.

\bibitem[R67]{rickert} N. W. Rickert, \emph{Amenable groups and groups with the fixed point property},  Trans. Amer. Math. Soc., Vol. 127, pp. 221-232, 1967. 

\bibitem[MR15]{matuirordam} H. Matui and M. Rørdam, \emph{Universal properties of group actions on locally compact spaces}, Journal of Functional Analysis, Vol. 268, No. 12, pp. 3601-3648, 2015.

\bibitem[M91]{banachlattice} Peter Meyer-Nieberg, \emph{Banach lattices}, Springer, Berlin, 1991.

\bibitem[M17]{monod} Nicolas Monod, \emph{Fixed points in convex cones}, Trans. AMS, B4, pp. 68-93, 2017.

\bibitem[N11]{navas} Andrés Navas,  \emph{Groups of Circle Diffeomorphisms}, Chicago Lecture in Mathematics Series, 2011.

\bibitem[P78]{palmer} T.W. Palmer, \emph{Classes of Nonabelian, Noncompact, locally compact groups}, Rocky Mountain Journal of Math., Vol.8, No.4, Fall, 1978.

\bibitem[P75]{plantefoliations}  J. F. Plante. \emph{Foliations With Measure Preserving Holonomy}, Annals of Mathematics - Second Series, Vol. 102, No. 2, pp. 327-361, September 1975.

\bibitem[P83]{plante} J.F. Plante, \emph{Solvable groups action on the line}, Trans. AMS, Vol. 274, No.1, July 1983.

\bibitem[P88]{paterson} Alan Paterson, \emph{Amenability}, Monograph AMS, Vol.29, 1988.

\bibitem[R72]{rosenblatt} Joseph Max Rosenblatt, \emph{Invariant linear functionals and counting conditions}, PhD thesis, University of Washington, Seattle, 1972.

\bibitem[R86]{rudinanal} Walter Rudin, \emph{Real and complexes analysis}, Higher Mathematics Series, McGraw-Hill Education, 3rd edition, May 1986.

\bibitem[R02]{runde} V. Runde, \emph{Lectures in Amenability}, Lecture Notes in Mathematics, Vol. 1774, Springer, Berlin Heidelberg, 2002.

\bibitem[S74]{banachlatticeandposoperator} H.H. Schaefer, \emph{Banach lattices and positive operators}, Die Grundlehren der mathematischen Wissenschaften, Vol. 215, Springer, 1974.

\end{thebibliography}
\end{document}